\documentclass[
    crop=false,
    a4paper,
    10pt,
]{article}

\usepackage[utf8]{inputenc} 
\usepackage[T1]{fontenc}
\usepackage{custom-style}
\usepackage{math}
\usepackage[plain]{theoremtemplate-style}
\usepackage[margin=2cm]{geometry}

\title{Multisymplectic observable reduction using constraint triples} 
\author{
	Antonio Michele MITI 
	\thanks{Dipartimento di Matematica, Sapienza Università di Roma, Piazzale Aldo Moro 5, 00185 Roma, Italy.
	\href{mailto:antoniomichele.miti@uniroma1.it}{antoniomichele.miti@uniroma1.it}}
	~ and Leonid RYVKIN
	\thanks{Institut Camille Jordan, Université Claude Bernard Lyon 1, 
			43 boulevard du 11 novembre 1918, 69622 Villeurbanne France
	\href{mailto:leonid.ryvkin@math.univ-lyon1.fr}{leonid.ryvkin@math.univ-lyon1.fr}.}
}

\begin{document}
\maketitle

\begin{abstract}
    The purpose of this paper is to present a fully algebraic formalism for the construction and reduction of $L_\infty$-algebras of observables inspired by multisymplectic geometry, using Gerstenhaber algebras, BV-modules, and the constraint triple formalism.  
    In the "geometric case", we reconstruct and conceptually explain the recent results of \cite{blacker2023reduction}.
    %
\end{abstract}

\tableofcontents

\section*{Introduction}\label{sec:intro}
\addcontentsline{toc}{section}{\nameref{sec:intro}}

Reduction can be described, in very broad strokes, as a procedure that associates to a given object (such as a manifold or an algebra) a “smaller” object of the same kind.
The main mechanisms by which this can happen are by looking at \emph{subspaces} or by dividing out symmetries, that is, taking \emph{quotients}. 
In the context of mechanical systems admitting some version of the Noether theorem, i.e. that have a link between symmetries and conserved quantities, both mechanisms are usually applied at the same time, in the form of a \emph{subobject-quotient} reduction scheme.

One of these settings that admits such reductions by subobject-quotients is symplectic geometry, and its generalisation, multisymplectic geometry. In the latter, a geometric reduction procedure was presented in \cite{blackerReductionMultisymplecticManifolds2021} and an algebraic one in \cite{blacker2023reduction}.\\ 
In this article, we provide a conceptual explanation for the results of  \cite{blacker2023reduction} by interpreting them using the constraint triples originally introduced by Dippell, Esposito and Waldmann \cite{dippelespositowaldmann} as an algebraic framework for coisotropic reduction. In doing so, our constructions become very general: We obtain a procedure for creating (and reducing) an $L_\infty$-algebra (of observables) given any Gerstenhaber algebra equipped with a BV-module and a cocycle.

\paragraph{Geometric reduction scheme.}
	The starting point of this work is the celebrated Marsden--Weinstein--Meyer  (MWM) reduction theorem~\cite{Meyer73,MarsdenWeinstein74}, a cornerstone of symplectic geometry providing a procedure to dimensionally reduce the data encoded by a given symplectic manifold in the presence of symmetry. 
	More technically, the theorem applies in the following setting:
\begin{itemize}
    \item A smooth manifold $M$ equipped with a symplectic form $\omega \in \Omega^2(M)$.
    \item A smooth action $\vartheta: G \times M \to M$ of a Lie group $G$ preserving the symplectic structure, that is, $\vartheta_g^* \omega = \omega$ for all $g \in G$.
    \item An equivariant momentum map $\mu: \g \to C^\infty(M)$ satisfying:
    \begin{itemize}
        \item $\d(\mu(\xi)) = \iota_{\underline{\xi}} \omega$ for all $\xi \in \g$, where $\underline{\xi}$ denotes the fundamental vector field generated by $\xi$, given by $\underline{\xi}(m) = \left.\frac{d}{dt}\right|_{t=0} \vartheta_{\exp(t\xi)}(m)$.
        \item $\mu(\mathrm{Ad}_{g^{-1}}\xi) = \mu(\xi) \circ \vartheta_g$ for all $g \in G$, where $\mathrm{Ad}$ denotes the adjoint representation of $G$ on its Lie algebra $\g$.
    \end{itemize}
    \item The group action is free and proper on the zero level set $\hat{\mu}^{-1}(0)$, where the associated momentum map $\hat{\mu}: M \to \g^*$ is defined by $\hat{\mu}(m)(\xi) = \mu(\xi)(m)$. This condition implies that $0 \in \g^*$ is a regular value of $\hat{\mu}$.
\end{itemize}
	Under these hypotheses, the MWM theorem asserts that the quotient $M_{\red} = \hat{\mu}^{-1}(0)\slash G$ is a smooth manifold that inherits a uniquely determined symplectic form $\omega_{\red}$ such that $\pi^* \omega_{\red} = i^* \omega$, where $\pi: \hat{\mu}^{-1}(0) \to M_{\red}$ is the canonical projection and $i: \hat{\mu}^{-1}(0) \hookrightarrow M$ is the inclusion map.
	Such a result has several far-reaching consequences. 
	It provides, for instance, the symplectic structure on coadjoint orbits, realised as reductions of the cotangent bundle $T^*G$ under the conjugation action of $G$ on itself. 
	More generally, it underlies cotangent bundle reduction: given a configuration space $Q$ with a free and proper $G$-action, the canonical Hamiltonian $G$-action on the phase space $M=T^*Q$ admits a momentum map, and the reduced space at zero momentum is symplectomorphic to $T^*(Q\slash G)$. 
	In other words, the theorem provides a geometric formulation of Noether’s principle:  
    Whenever a Hamiltonian system with symmetry group $G$ admits a momentum map, the conserved quantities associated with the symmetries are encoded in its image. Moreover, the dynamics of the system descends to the reduced phase space $M \slash G$, obtained via symplectic reduction.

\paragraph{Limitations of the classical reduction scheme.}
	The above setting relies on strong structural assumptions, in particular the freeness and properness of the group action, as well as the regularity of the momentum level set. 
	These conditions often fail in applications.
	\\
	To address these limitations, various algebraic approaches to reduction have been proposed in the literature. These methods shift the focus from manifolds to the algebra of observables, aiming to avoid direct quotient constructions by instead encoding symmetry and constraint data at the algebraic level. 
	This perspective was already present in the early development of the field (see~\cite{SniatWein83} and~\cite{MarsdenWeinstein01,Sniatycki2005} for modern expositions), and has recently been reformulated in terms of the so-called \emph{constraint triples}.

\paragraph{Constraint triples and the algebraic approach to reduction.}
	The theory of \emph{constraint triples} has been developed in recent publications~\cite{dippelespositowaldmann,Dippell2023,dippelkern} (and references therein). 
	The guiding idea stems from coisotropic reduction and aims to avoid the potential pitfalls of computing quotients directly by instead encoding all reduction data into a single categorical object that captures the full structure involved.
	\\
	The key observation is that a general reduction scheme can be formalized as a \emph{subobject–quotient} process: starting from a total space $M_\Total$, one first restricts to a subspace $M_\Wobs \subset M_\Total$, and then imposes an equivalence relation $\sim_\Null$ on $M_\Wobs$, obtaining the reduced space $M_\red:= M_\Wobs/\sim_\Null$.
	For instance, in classical MWM reduction, the total space is a symplectic manifold $M_\Total = M$, the subspace is the zero level set $M_\Wobs = \mu^{-1}(0)$ of the moment map, and the equivalence relation is induced by the Lie group action, i.e.\ $\sim_\Null = \sim_G$.
	\\
	A naive attempt to reproduce this reduction purely at the level of observables (i.e., smooth functions) may proceed as follows:
	\begin{itemize}
		\item First restrict to $M_\Wobs$ by passing to the quotient $C^\infty(M_\Total) \to C^\infty(M_\Total)/I_{M_\Wobs}$, where $I_{M_\Wobs} := \{ f \in C^\infty(M_\Total) \mid f|_{M_\Wobs} = 0 \}$ is the vanishing ideal of $M_\Wobs$;
		\item Then impose the group action by taking $G$-invariants: $\left(C^\infty(M_\Total)/I_{M_\Wobs}\right)^G$.
	\end{itemize}
	Note that this algebraic procedure follows a \emph{quotient–subobject} pattern, in contrast to the \emph{subobject–quotient} structure of the MWM theorem. However, the two approaches can be reconciled by rewriting the invariant quotient as
	\begin{align}\label{eq:algredexchange}
		\left(\frac{C^{\infty}(M_\Total)}{I_{M_\Wobs}}\right)^G 
		= \frac{\{ f \in C^{\infty}(M_\Total) \mid \forall g \in G: \vartheta_g^*f - f \in I_{M_\Wobs} \}}{I_{M_\Wobs}},
	\end{align}
	where $\vartheta_g^*$ denotes the pullback along the action of $g \in G$.  
	This formulation reveals the subquotient nature of the construction also at the algebraic level (cf.~\autoref{eq:algredexchange} and the singular reduction schemes of~\cite{Sniatycki2005}).

\paragraph{Multisymplectic reduction schemes}
	Historically, the transition from symplectic to multisymplectic geometry was motivated by the limitations of applying the symplectic formalism to field theory. 
	While one can formally equip the (infinite-dimensional) space of solutions to the field equations with a presymplectic structure, this approach faces substantial technical obstacles, starting from the difficulties in pinpointing the smooth structure on the ensuing infinite-dimensional submanifolds.
	The covariant multisymplectic formalism avoids these difficulties by working instead with a finite-dimensional \emph{multiphase space}, a multisymplectic manifold naturally associated to the configuration bundle (We refer to ~\cite{Gimmsy1,Blohmann2021,Ryvkin2018} for background on the geometric formulation of classical field theory that underlies this framework).
	\\
	For the sake of simplifying the treatment of infinite-dimensional mechanical systems and field theories with symmetries, the problem of extending reduction theorems to the multisymplectic context has garnered significant interest.
	A reduction scheme for multisymplectic Hamiltonian $G$-spaces is proposed in~\cite{blackerReductionMultisymplecticManifolds2021}, following the spirit of the Marsden--Weinstein--Meyer theorem and thus restricted to a relatively regular setting. 
	More recently, in~\cite{blacker2023reduction}, we introduce a purely algebraic version of the latter wich is also applicable to singular cases. This approach is inspired by the formal adaptation of classical symplectic reduction techniques to the $L_\infty$-algebra of observables associated with a multisymplectic manifold~\cite{rogers2012linfty}, and operates directly at the level of observables. It accommodates arbitrary group actions and covariant momentum maps, independently of regularity assumptions.
	In the present work, we further develop this algebraic perspective by introducing a multisymplectic reduction scheme based on the formalism of constraint triples~\cite{dippelespositowaldmann}. The latter provides a more conceptual, less ad-hoc derivation of the observable reduction scheme in~\cite {blacker2023reduction}.

\subsection*{Main Results and Outlook}
%
Guided by the aim of providing a deeper understanding of the algebraic reduction scheme proposed in~\cite{blacker2023reduction}, we revisit the algebraic structures underlying Cartan calculus and, by extension, the multisymplectic observable algebra introduced by Rogers in~\cite{rogers2012linfty}.
In \autoref{sec:BV-Cartan-calculus}, borrowing the formalism of Tsygan–Tamarkin–Nest noncommutative differential calculus, we introduce a purely algebraic formulation of the multisymplectic observable algebra. Namely, we show that, when a certain cochain is chosen, any BV-module gives rise to an associated $L_\infty$-algebra of observables.
\begin{bigthm}[\autoref{sec:linftyfrombv}]\label{bigthm:AlgebraicRogersAlgebra}
	Let \( (G, \wedge, \{ \blank, \blank \}) \) be a Gerstenhaber algebra, and \( (V, \d) \) a BV-module over \( G \). Let \( \omega \in V^{k+1} \) be a cocycle, i.e., \( \d\omega = 0 \), with $k\geq 1$.
	\\
	The graded vector space $\Ham(V, \omega) := \bigoplus_{i \leq 0} \Ham^i(V, \omega)$ where
	$$
	\Ham^0(V, \omega) := \left\{ (\alpha, X) \in V^{k-1} \oplus G^1 \mid \iota_X \omega = -\d\alpha \right\},
	\qquad
	\Ham^i(V, \omega) := V^{k-1+i} \quad \text{for } i < 0~,
	$$
	admits an \( L_\infty \)-algebra structure with multibrackets $ l_j \colon \Ham^{\otimes j}(V, \omega) \to \Ham(V, \omega)[2-j]$ given by \autoref{def:AlgebraicRogersAlgebra}.
\end{bigthm}
In \autoref{sec:constraint-Cartan-Calculus}, we develop the constraint algebra framework, introducing the notion of constraint BV-modules (\autoref{def:constraintbv}), showing how the above-mentioned noncommutative differential calculus framework naturally integrates with the algebraic theory of constraint triple. 
In particular, we prove in \autoref{lem:const-bv-from-const-LR} how such structures could naturally arise from constraint Lie–Rinehart algebras, thereby generalising the constraint Cartan calculus introduced in~\cite{Dippell2023,dippelkern}.
\begin{bigthm}[\autoref{lem:const-bv-from-const-LR}]\label{bigthm:const-bv-from-const-LR}
	Let $\mathfrak{L}$ be a constraint Lie--Rinehart algebra over a constraint algebra $A$ as per \cite[Def. 4.1]{dippelkern}. 
	Then the constraint vector space $\Lambda_A \mathfrak L$ is constraint Gerstenhaber algebra and $\CE(\mathfrak L) = ((\Lambda_A \mathfrak{L})^*)^\rev$ is a constraint BV-module over it.
\end{bigthm}
In the final section, we unify the previous two constructions. 
In \autoref{thm:mainredgenalgleotired}, we show how to construct a constraint $L_\infty$-algebra of observables from any constraint BV-module. 
\begin{bigthm}[\autoref{lem:clinftyconstr}]\label{bigthm:clinftyconstr}
	Let $(G, \wedge, \{ \blank, \blank \})$ be a constraint Gerstenhaber algebra, and $(V, \d)$ a constraint BV-module over $G$.
	Fix a $k+1$-cocycle $\omega \in V^{k+1}_\Wobs$ with $\d\omega = 0$.
	The \emph{constraint graded vector space}
	\[
		\Ham(V, \omega) := \bigoplus_{i \leq 0} \Ham^i(V, \omega),
	\]
	given in degree $0$ by the constraint vector space of Hamiltonian pairs \( \Ham^0(V, \omega) \) as per \autoref{def:constraint-algebraicham0}, and by $\Ham^i(V, \omega) := V^{k-1 + i}$ for any $i \leq -1$, forms a constraint $L_\infty$-algebra when endowed with the family of constraint multilinear maps \( l_j \colon \Ham^{\otimes j}(V, \omega) \to \Ham(V, \omega)[2-j] \) defined as in \autoref{def:AlgebraicRogersAlgebra}.
\end{bigthm}
This general result is then applied to specific cases. 
In particular, in \autoref{thm:mainredgenalgleotired} and \autoref{cor:howtorecoverourreduction} we recover the reduction scheme of~\cite{blacker2023reduction} as a special instance of our framework and clarify the role of the so-called “residue defect” within this algebraic setting in \autoref{ssec:geoRed-residue-defect}.

\paragraph{Outlooks.}
	In recent years, several reduction formalisms for classical field theory have been discussed that resonate with the present framework based on constraint triples.
	Without any claim of completeness, we mention:
	\begin{itemize}
		\item In the context of \emph{polysymplectic} geometry, which generalizes multisymplectic geometry by considering closed vector-valued $2$-forms, the MWM-reduction scheme has been adapted in \cite{Marrero2015} Marrero-Román--Roy-Salgado-Vilariño. Recently, De Lucas–Rivas–Vilariño–Zawora \cite{DeLucas2023} refined the previous work and discussed application to \emph{k-polycosymplectic} Hamiltonian systems with field symmetries, and Capriotti-Díaz-García-Toraño Andrés-Mestdag \cite{Capriotti2022} discussed Lagrangian field theories and the reconstruction process.
		\item  In \cite{Bua2015}, Búa–Mestdag–Salgado investigated the reduction process of a \emph{$k$-symplectic} (i.e. a family of $k$ closed two-forms with a suitably compatible distribution) field theory whose Lagrangian is invariant under a symmetry group, showing that under suitable regularity conditions one can define a corresponding mechanical connection.
		\item Very recently, De Lucas–Rivas–Vilariño–Zawora \cite{DeLucas2025} introduced a MWM-type reduction for \emph{$k$-contact} field theories, extending reduction schemes to dissipative systems.
		\item On the applied side of the multisymplectic approach to infinite dimensional mechanical systems, de Lucas-Gràcia-Rivas-Román--Roy-Vilariño \cite{DeLucas2022} emploied higher degree closed forms to the study of the reduction of \emph{Lie systems} (a particular class of non-autonomous PDEs). Moreover, Gomis–Guerra–Román--Roy \cite{Gomis2023a} developed the multisymplectic formalism for field theories in relevant physical models, focusing on reduction in the presence of mechanical constraints rather than pure symmetries.
\end{itemize}
A natural evolution of the present work is to clarify the relation between these approaches and the constraint $L_\infty$-algebraic formalism developed here and in~\cite{blacker2023reduction}.
Notably, the $L_\infty$-algebra of observables was used in Bernardy’s recent work~\cite{Bernardy2025} for reduction in the context of higher‐order Lagrangian field theory \cite{Blohmann2021,Blohmann24}. It would be interesting to explore whether there is a connection between the reduction (which is phrased in terms of homotopy momentum maps, originally introduced in \cite{CalliesFregierRogersZambon16}) and ours.

More generally, our work suggests that the constraint formalism provides a powerful tool to treat singular geometric (and algebraic) settings. Rather than defining all algebraic structures in the setting of constraint triples "by hand", it could be interesting to provide a systematic treatment, e.g. within the language of operads.

\subsection*{Structure of the paper}
The manuscript is organised into four sections.

\autoref{subsec:multisymIntro} is a self-contained introduction to multisymplectic geometry, providing the necessary background on Hamiltonian pairs, the observable $L_\infty$-algebra and covariant momentum maps.

\autoref{sec:BV-Cartan-calculus} contains a short survey on the algebraic framework of Cartan calculus in terms of Batalin–Vilkovisky (BV) modules, and provides a recipe to associate an observable $L_\infty$-algebra with any BV-module with a fixed cocycle (see \autoref{bigthm:AlgebraicRogersAlgebra}).

In \autoref{sec:constraint-Cartan-Calculus} we recall the formalism of constraint algebras and modules, introducing the notion of constraint BV-modules and showing how they naturally arise from constraint Lie–Rinehart algebras (see \autoref{bigthm:const-bv-from-const-LR}).

Finally, \autoref{sec:multisymplecticconstraintreduction} introduces the notion of constraint $L_\infty$-algebra of observables (see \autoref{bigthm:clinftyconstr}) and defines their reduction, before establishing the comparison with the geometric case.

\subsection*{Notations and conventions}
%
We work over a fixed ground field $\mathbb{k}$ of characteristic zero, and consider the category of $\mathbb{Z}$-graded vector spaces over $\mathbb{k}$ seen as the functor category from the discrete category $\mathbb{Z}$ to the category of vector spaces over $\mathbb{k}$.
\\
Whenever no confusion arises, we identify a graded vector space $V$ with the direct sum of its homogeneous components:
\[
V = \bigoplus_{i \in \mathbb{Z}} V^i.
\]

Given a graded vector space $V$, we denote by $V^{\rev}$ its \emph{reverse grading}, defined by
\[
V^{\rev} := \bigoplus_{i \in \mathbb{Z}} V^{-i}.
\]

All differential graded vector spaces are understood to follow the \emph{cohomological convention}, i.e., they are cochain complexes with a differential $\d$ of degree $+1$:
\[
\d \colon V^i \longrightarrow V^{i+1}, \quad \d^2 = 0.
\]

Throughout this work, we consider algebraic structures defined on graded vector spaces. We will typically omit the qualifier $\mathbb{Z}$-graded, as most vector spaces under consideration are graded by default; ungraded vector spaces are implicitly regarded as graded and concentrated in degree zero. 
Similarly, all commutative algebras are assumed to be associative, and all associative algebras to be unital. These qualifiers will generally be omitted in the text for the sake of readability.

Unless otherwise specified, tensor products $\otimes$ are taken over the base field~$\mathbb{k}$. When dealing with $A$-modules, tensor products over the algebra~$A$ may also appear; in such cases, we will explicitly write $\otimes_A$.

\subsection*{Acknowledgements}
	The authors thank Casey Blacker, Peter Crooks, Marco Valerio D’Agostino, Marvin Dippell, Chiara Esposito, Domenico Fiorenza, David Kern, Stefan Waldmann, and Luca Vitagliano for their many fruitful discussions, which contributed significantly to the development of this manuscript. The authors also thank the editor and referee for their careful reading of the manuscript.
	LR was funded by the Deutsche Forschungsgemeinschaft (DFG, German Research
	Foundation) – Projektnummer 539126009
	and acknowledges the financial support of the LabEx MiLyon. He thanks Sapienza University of Rome for its hospitality during his research stay.\\
	AMM gratefully acknowledges funding from the European Union’s \emph{Horizon 2020 Research and Innovation Programme} under Grant Agreement No.~\emph{101034324}, and partial support from the \emph{Italian Group for Algebraic and Geometric Structures and their Applications} (GNSAGA–INdAM).  
	The author also thanks the \emph{Dipartimento di Matematica dell’Università di Salerno} and the \emph{Institut Camille Jordan} of \emph{Université Lyon 1} for their hospitality during the preparation of this work.
	\\  
	In accordance with Open Access policies, a \emph{CC BY-NC-SA} public copyright license has been applied by the authors to this document. It will extend to all subsequent versions, including the \emph{Author Accepted Manuscript} resulting from this submission.

\section{Overview on multisymplectic manifolds and Hamiltonian pairs}\label{subsec:multisymIntro}
%
Multisymplectic geometry generalises the concepts of symplectic geometry by considering closed, nondegenerate differential forms of degree greater than two, providing a natural framework for field theories that extends the Hamiltonian formalism of classical mechanics.
In this section, we review the notion of multisymplectic manifolds, along with their associated observable algebras and momentum maps.

\begin{definition}[(Pre-)multisymplectic manifolds]
	Let $k \geq 1$.
	A \emph{premultisymplectic manifold} of degree $k$ is a smooth manifold $M$ endowed with a closed differential form $\omega \in \Omega^{k+1}(M)$.
	The manifold is said to be \emph{multisymplectic} (or \emph{$k$-plectic}) if, in addition, $\omega$ is non-degenerate in the sense that the bundle map
	\[
		\iota_\blank \omega: TM \to \Lambda^k T^*M,\quad v \mapsto \iota_v \omega
	\]
	is injective.
\end{definition}

\begin{example}
	We list several classes of multisymplectic manifolds:
	\begin{itemize}
		\item A \emph{symplectic manifold} $(M, \omega)$, where $\omega \in \Omega^2(M)$ is a closed non-degenerate 2-form, is a multisymplectic manifold of degree $k=1$ .

		\item A smooth manifold $M$ of dimension $n \geq 1$ endowed with a \emph{volume form} $\eta \in \Omega^n(M)$ is a multisymplectic manifold of degree $k = n-1$.
		      Non-degeneracy follows from the fact that a volume form is nowhere vanishing.

		\item Let $Q$ be a smooth manifold, and consider the manifold $M = \Lambda^k T^*Q$, i.e., the bundle of $k$-forms over $Q$.
		      There exists a canonical $k$-form $\theta$ on $M$ defined by
		      \[
			      \theta_\alpha(v_1,\dots,v_k) = \alpha(T\pi(v_1),\dots,T\pi(v_k)),
		      \]
		      where $\pi: M \to Q$ is the canonical projection. The form $\omega = \d\theta$ endows $M$ with a multisymplectic structure of degree $k$. This construction plays a central role in covariant Hamiltonian formulations of classical field theory.

		\item The trivial ($0$) $k$-form on any manifold is premultisymplectic but not multisymplectic, as the non-degeneracy condition fails.
	\end{itemize}
\end{example}

On a symplectic manifold $(M, \omega)$, every smooth function $f \in C^\infty(M)$ uniquely determines a vector field $X_f$ via the \emph{Hamiltonian relation} $\iota_{X_f} \omega = -\d f$.
A natural attempt to generalize this correspondence to the multisymplectic setting is to replace the function $f$ with a differential form $\alpha \in \Omega^{k-1}(M)$ and seek a vector field $X$ such that
\[
	\d \alpha = -\iota_X \omega~.
\]
However, in contrast to the symplectic case, neither the existence nor the uniqueness of such a vector field is guaranteed in general. In particular, on a premultisymplectic manifold, a given form $\alpha$ may correspond to multiple vector fields, or none at all. Consequently, it becomes necessary to explicitly consider pairs $(\alpha, X)$, recording both the form and the associated vector field. This observation motivates the introduction of the following algebraic structure.

\begin{definition}[Hamiltonian Leibniz algebra {\cite{rogers2012linfty}}]\label{def:ham0}
	Let $(M, \omega)$ be a premultisymplectic manifold of degree $k$. We define:
	\begin{enumerate}
		\item The \emph{space of Hamiltonian pairs} as
		      \[
			      \Ham^0(M, \omega) := \left\{ (\alpha, X) \in \Omega^{k-1}(M) \times \X(M) \;\big|\; d\alpha = -\iota_X\omega \right\}.
		      \]
		\item The \emph{Leibniz bracket} $[\blank, \blank]: \Ham^0(M, \omega)^{\otimes 2} \to \Ham^0(M, \omega)$ of Hamiltonian pairs by
		      \[
			      [(\alpha_1, X_1), (\alpha_2, X_2)] = (\Lie_{X_1} \alpha_2, [X_1, X_2])~,
			      \qquad\qquad \forall (\alpha_i,X_i) \in \Ham^0(M,\omega)~.
		      \]
	\end{enumerate}
	The pair $(\Ham^0(M, \omega), [\blank, \blank])$ forms a left Leibniz algebra:  
	\[
		[a,[b,c]] = [[a,b],c] + [b,[a,c]]~,
		\qquad\qquad \forall a,b,c \in \Ham^0(M,\omega)~.
	\]
\end{definition}
For \(k > 1\), the Leibniz bracket \([\blank, \blank]\) is generally not skew-symmetric and therefore does not define a Lie algebra. Moreover, a naive skew-symmetrisation of the Leibniz bracket on Hamiltonian pairs typically fails to satisfy the Jacobi identity.
In contrast, for \(k = 1\), one recovers the classical Lie algebra of Hamiltonian vector fields.

Nevertheless, one can define a natural skew-symmetric bracket on Hamiltonian pairs:
\begin{align*}
	\{(\alpha, X), (\beta, Y)\} := (\iota_Y \iota_X \omega, [X, Y]).
\end{align*}
For $k > 1$, this bracket does not satisfy the Jacobi identity; instead, its Jacobiator is exact. In fact, this bracket can be extended to an $L_\infty$-algebra, first constructed in~\cite{rogers2012linfty}.

\begin{definition}[Observables $L_\infty$-algebra {\cite[Sec.~4]{rogers2012linfty}}]\label{def:linftygeom}
	Let $(M, \omega)$ be a pre-multisymplectic manifold of degree $k$, i.e., \( \omega \in \Omega^{k+1}(M) \) with \( \d\omega = 0 \).
	We define the \emph{graded space of multisymplectic observables} as
	\[
		\Ham(M, \omega) := \bigoplus_{i = -k+1}^0 \Ham^i(M, \omega),
	\]
	where
	\[
		\Ham^0(M, \omega) := \left\{ (\alpha, X) \in \Omega^{k-1}(M) \oplus \X(M) \mid \iota_X \omega = -\d\alpha \right\}
	\]
	is the space of Hamiltonian pairs (\autoref{def:ham0}), and for each \( i \in \{-k+1, \dots, -1\} \),
	\[
		\Ham^i(M, \omega) := \Omega^{k-1+i}(M).
	\]
	\noindent
	This graded vector space carries an $L_\infty$-algebra structure with degree $(2-j)$ multibrackets $l_j\colon \Ham(M,\omega)^{\otimes j} \to \Ham(M,\omega) [2 - j]$, defined as:
	\begin{itemize}
		\item \( l_1 \colon \Ham^i(M,\omega) \to \Ham^{i+1}(M,\omega) \):
		      \[
			      l_1(\alpha) = \begin{cases}
				      \d\alpha      & \text{if } i < -1~, \\
				      (\d\alpha, 0) & \text{if } i = -1~, \\
				      0             & \text{if } i = 0~;
			      \end{cases}
		      \]
		\item \( l_2 \colon \Ham^0(M,\omega) \otimes \Ham^0(M,\omega) \to \Ham^{0}(M,\omega) \):
		      \[
			      l_2\big( (\alpha, X), (\beta, Y) \big) = (\iota_X \iota_Y \omega,~ [X, Y])~;
		      \]
		\item For each \( j \in \{3, \dots, k+1\} \):
		      \[
			      l_j\big( (\alpha_1, X_1), \dots, (\alpha_j, X_j) \big) = -\iota_{X_1} \dots \iota_{X_j} \omega~;
		      \]
		\item All multibrackets of arity greater than or equal to $2$ vanish when evaluated on inputs containing at least one element of nonzero degree.
	\end{itemize}
	We refer to $ \Ham(M, \omega) $ endowed with these brackets as the \emph{$L_\infty$-algebra of multisymplectic observables}\footnote{Observe that here we are using a different convention than Rogers, namely we are contracting sequences of vector fields in the reverse order. The $(-1)^{\frac{j(j+1)}{2}}$ that appears in \cite{rogers2012linfty} is the Koszul sign of the permutation reversing the order of a sequence of degree $1$ elements $X_1,\dots, X_j$.}
    .
\end{definition}

Intuitively, we can think of an $L_\infty$-algebra as a higher analogue of a (graded) Lie algebra, where the Jacobi identity and its higher analogues hold up to homotopy (i.e. up to exact terms). Since the $L_\infty$-algebra in \autoref{def:linftygeom} is of a very special kind (most of its brackets are non-trivial only when restricted to $\Ham^0(M,\omega)$),  the $L_\infty$-algebra relations (generalised Jacobi identities) simplify. Explicitly, they are given by:
\begin{align*}
	\partial l_j = l_1 \circ l_{j+1} && 2\leq  j \leq k+1
\end{align*}
where we set \( l_{k+2} := 0 \) and \( \partial \) denotes the Chevalley–Eilenberg coboundary:
\[
	(\partial l_m)(x_1, \dots, x_{m+1}) = \sum_{1 \leq i < j \leq m+1} (-1)^{i+j} l_m\bigl(l_2(x_i, x_j), x_1, \dots, \widehat{x_i}, \dots, \widehat{x_j}, \dots, x_{m+1}\bigr).
\]
These identities are consequences of \( \d\omega = 0 \) and standard Cartan calculus.

\begin{remark}
	If $\omega$ is non-degenerate, the vector field component of a Hamiltonian pair is uniquely determined by its form component.
	Even in the degenerate case, an $L_\infty$-structure on forms alone can be obtained by choosing, for each form, a compatible vector field completing it to a Hamiltonian pair.
	The relationship between this construction and \autoref{def:linftygeom} can be found in \cite[\S 4]{CalliesFregierRogersZambon16}.
\end{remark}

\begin{remark}
	This construction is closely related to a general method of endowing free resolutions of Lie algebras with $L_\infty$-structures (cf.~\cite{barnichShLieStructure1998}).
	In the multisymplectic setting, the part of the de Rham complex resolving the space of Hamiltonian vector fields is usually only locally free (as a sheaf). However, the construction from \cite{barnichShLieStructure1998} still works and leads to the $L_\infty$-algebra $\Ham(M,\omega)$. 
	In the special case of volume forms, this resolution was already presented in~\cite{rogerExtensionsCentralesAlgebres1995} for the study of central extensions.
\end{remark}

\begin{example}
	Consider \( M = \mathbb{R}^3 \) with \( \omega \) the standard volume form.
	The associated $L_\infty$-algebra is
	\[
		C^\infty(\mathbb{R}^3) \oplus \mathrm{Ham}^0(M, \omega),
	\]
	where \( \mathrm{Ham}^0(M, \omega) \subset \Omega^1(\mathbb{R}^3) \times \mathfrak{X}(\mathbb{R}^3) \).
	Since \( \omega \) is non-degenerate, the 1-form uniquely determines the vector field component, and
	\( \Lambda^2 T^*\mathbb{R}^3 \cong T\mathbb{R}^3 \) ensures that any 1-form defines a vector field.
	Thus, the underlying space reduces to
	\[
		C^\infty(\mathbb{R}^3) \oplus \mathfrak{X}(\mathbb{R}^3).
	\]

	Identifying the de Rham complex as
	\[
		C^\infty(\mathbb{R}^3) \xrightarrow{\operatorname{grad}} \mathfrak{X}(\mathbb{R}^3) \xrightarrow{\operatorname{rot}} \mathfrak{X}(\mathbb{R}^3) \xrightarrow{\operatorname{div}} C^\infty(\mathbb{R}^3),
	\]
	we describe the multisymplectic $L_\infty$-algebra on \( L^{-1} \oplus L^0 = C^\infty(\mathbb{R}^3) \oplus \mathfrak{X}(\mathbb{R}^3) \) as:
	\begin{itemize}
		\item \( l_1(f) = \operatorname{grad}(f) \);
		\item \( l_2(X, Y) = \operatorname{rot}(X) \times \operatorname{rot}(Y) \);
		\item \( l_3(X, Y, Z) = -\bigl( (\operatorname{rot}(X) \times \operatorname{rot}(Y)) \cdot \operatorname{rot}(Z) \bigr) \).
	\end{itemize}

	These satisfy:
	\begin{itemize}
		\item \( l_2(l_1(f), Y) = 0 \);
		\item \( l_2(l_2(X, Y), Z) + \text{cyclic} = l_1 l_3(X, Y, Z) \);
		\item \( l_3(l_2(X, Y), Z, W) + \text{cyclic} = 0 \).
	\end{itemize}

	The Hamiltonian vector field of \( X \in L^0 \) is \( \operatorname{rot}(X) \), which intertwines \( l_2 \) with the Lie bracket.
	Higher brackets depend only on the \( \operatorname{rot} \) of the inputs.
	Notice that working directly on the de Rham side simplifies the verification of the $L_\infty$ relations.
	(In contrast, computations in terms of \( \operatorname{grad} \), \( \operatorname{rot} \), and \( \operatorname{div} \) are long and unintuitive.)
\end{example}

\subsection{Symmetries and reduction}\label{ssec:symred}
%
In the context of multisymplectic geometry, the notion of symmetry—traditionally understood in symplectic geometry as a Lie group action preserving the symplectic form—admits a natural generalisation to group actions that preserve the multisymplectic form. 
Based on the Leibniz algebra structure previously introduced, one may accordingly define momentum maps in the premultisymplectic setting as follows.
\begin{definition}[Covariant momentum map, {\cite[Def.~3.2]{blackerReductionMultisymplecticManifolds2021}}]
	\label{def:covmom}
	Let $(M, \omega)$ be a premultisymplectic manifold, and let $\vartheta \colon G \times M \to M$ be a $G$-action preserving $\omega$.
	Denote by $\underline{\blank} \colon \g \to \X(M)$ the induced infinitesimal action.
	A \emph{covariant momentum map} is a morphism of Leibniz algebras
	\begin{displaymath}
		\mu \colon \g \to \Ham^0(M, \omega)
	\end{displaymath}
	satisfying $\pi_{\X(M)}(\mu(\xi)) = \underline{\xi}$ for all $\xi \in \g$, where $\pi_{\X(M)}$ is the projection $\Ham^0(M, \omega)\subset\Omega^{k-1}(M) \times \X(M) \to \mathfrak{X}(M)$.
\end{definition}
\noindent 
We also write $\pi_{\Omega(M)}$ for the projection $\Ham^0(M, \omega)\subset\Omega^{k-1}(M) \times \X(M) \to \Omega^{k-1}(M)$. Moreover, we will, by abuse of notation, identify $\mu$ with its composition $\pi_{\Omega(M)} \circ \mu$, since the infinitesimal action always determines the vector field component.

The notion of momentum maps plays a central role in the formulation of reduction procedures, which produce lower-dimensional manifolds that are equivalent, in a suitable sense, to the original ones under the action of a symmetry group. 
In the multisymplectic setting, a geometric reduction scheme has been developed in~\cite{blackerReductionMultisymplecticManifolds2021} and it requires the following data:
\begin{itemize}
    \item A multisymplectic manifold \((M, \omega)\) of degree \(k\).
    \item A smooth action \(\vartheta: G \times M \to M\) of a Lie group \(G\) preserving the multisymplectic form \(\omega\).
    \item A covariant momentum map \(\mu\), which is \(G\)-invariant when regarded as a map \(\hat{\mu}: M \to \Lambda^k T^*M \boxtimes \mathfrak{g}^*\), where $\boxtimes$ denotes the external tensor product.
\end{itemize}
Theorem~1 in~\cite{blackerReductionMultisymplecticManifolds2021} asserts that if the zero level set \(\hat{\mu}^{-1}(0)\) is an embedded submanifold on which the \(G\)-action is free and proper, then the quotient \(\hat{\mu}^{-1}(0)/G\) inherits a natural premultisymplectic structure. It is worth noting that the result in \emph{loc.~cit.} is more general, allowing for reduction at arbitrary closed \(\mathfrak{g}^*\)-valued forms, provided the group action satisfies suitable regularity conditions.
In \autoref{sec:multisymplecticconstraintreduction} we discuss a reduction procedure that acts directly on the algebra $\Ham(M, \omega)$ and which is applicable to arbitrary Lie algebra actions and covariant momentum maps.

We conclude this subsection with the following two observations.
\begin{remark}
	Let $\mu$ be a covariant momentum map for a connected Lie group action. Then $\mu$ is automatically equivariant.
	In particular, in the symplectic case with connected group, the $G$-equivariance of $\hat\mu$ is equivalent to requiring that the map $f \mapsto X_f$ is a Lie algebra homomorphism, where $C^\infty(M)$ carries the Poisson bracket $\{f,g\} = \omega(X_f, X_g)$.
\end{remark}
\begin{remark} As we have seen, the Lie algebra structure of a symplectic manifold can be generalized in two different ways, either as an $L_\infty$-structure (i.e. upon relaxing the Jacobi identity - requiring it to hold only up to homotopy) or a Leibniz algebra structure (i.e. relaxing the skew-symmetry). Each of these choices has its associated notion of momentum map, a \emph{homotopy momentum map} as introduced in \cite{CalliesFregierRogersZambon16} for the $L_\infty$-perspective and the covariant momentum map as introduced above for the Leibniz perspective. In the reduction procedure described in \cite{blackerReductionMultisymplecticManifolds2021}, the Leibniz perspective is used, while the reduction approach of \cite{Bernardy2025} uses \emph{homotopy momentum maps}. Interestingly, for us, momentum maps of Leibniz type will be more useful, even though the objects we try to reduce are $L_\infty$-algebras.
\end{remark}


\section{An algebraic approach to the observables $L_\infty$-algebra}\label{sec:BV-Cartan-calculus}

In this section, we review the basic notions required to reformulate the classical Cartan calculus—concerning vector fields and differential forms on a smooth manifold—in a purely algebraic framework. 
Most of the material reviewed here is well-established, including its extensions to the noncommutative setting. The only novel contribution lies in \autoref{sec:linftyfrombv}, where we observe that the construction of the $L_\infty$-algebra introduced in~\cite{rogers2012linfty} extends to any BV-module structure equipped with a cocycle, and thus, in particular, to any Lie-Rinehart algebra.

\subsection{Cartan Calculus from an algebraic viewpoint}
A purely algebraic formulation of the Cartan calculus can be expressed through the language of BV-modules or Lie--Rinehart algebras.

Let ${A}$ be a commutative  algebra over a field $\mathbb{k}$.
Denote by $\X(A)=\Der_{\mathbb{k}}(A)$ its Lie algebra of derivations, which naturally forms an $A$-submodule of $\End(A)$ with the action $A\action \End(A)$ given by postmultiplying with the element of $A$.
\\
In the context of differential geometry, a primary example is when $A = C^\infty(M)$, the algebra of smooth functions on a smooth manifold $M$. In this case, $\mathfrak{X}(A)$ is isomorphic to the standard $C^\infty(M)$-module of smooth vector fields on $M$, often denoted by $\mathfrak{X}(M)$.

\subsubsection{Lie--Rinehart algebras}\label{subsec:LR}
Lie-Rinehart algebras were introduced as a natural generalisation of both Lie algebras and modules of derivations, allowing for a unified treatment of algebras of differential forms over arbitrary commutative rings.

\begin{definition}[Lie--Rinehart algebra \cite{rinehartDifferentialFormsGeneral1963}]
	Let $A$ be a commutative $\mathbb{k}$-algebra. A Lie--Rinehart algebra over $A$ is a pair $(\mathfrak{L}, \rho)$ where:
	\begin{itemize}
		\item $\mathfrak{L}$ is an $A$-module equipped with a $\mathbb{k}$-Lie algebra structure $[\blank,\blank]$,
		\item $\rho\colon \mathfrak{L} \to \X(A)$ is an $A$-linear Lie algebra homomorphism (called the anchor),
	\end{itemize}
	satisfying the Leibniz compatibility condition:
	\begin{align*}
		[X, fY] = f[X,Y] + \rho(X)(f)Y && \forall X,Y \in \mathfrak{L}, \, f \in A~.
	\end{align*}
	In the sequel, when working with different base algebras, we shall denote a Lie--Rinehart algebra simply by the pair $(A, \mathfrak{L})$, leaving the anchor map $\rho$ implicit unless explicitly needed.

	A morphism of Lie--Rinehart algebras $(A, \mathfrak{L}) \to (A', \mathfrak{L}')$ consists of an algebra morphism $\phi: A \to A'$ together with a Lie algebra morphism $\psi: \mathfrak{L} \to \mathfrak{L}'$ that is $A$-linear, where $\mathfrak{L}'$ is regarded as an $A$-module via $\phi$.

\end{definition}

\begin{example}\label{ex:XMisLR}
	The most basic example of a Lie--Rinehart algebra is the triple $(A, \X(A), \mathrm{id})$, where $\X(A)$ denotes the Lie algebra of $\mathbb{k}$-linear derivations of $A$, equipped with the identity anchor map.
    \\
	When $A = C^{\infty}(M)$ for a smooth manifold $M$, this recovers the classical Lie algebra of vector fields, and any involutive distribution $D \subset TM$ yields a Lie--Rinehart algebra structure on the $C^{\infty}(M)$-module of its sections $\Gamma(D)$.
    \\
	More generally, singular foliations in the sense of Androulidakis and Skandalis can be viewed as Lie--Rinehart algebras, as their space of vector fields is closed under the Lie bracket and stable under multiplication by smooth functions~\cite{androulidakisHolonomyGroupoidSingular2009}.
    \\
	Finally, a Lie--Rinehart algebra over $C^{\infty}(M)$ that is projective as a module is equivalent to a Lie algebroid, as established in \cite{Huebschmann1998}.
\end{example}

The anchor map $\rho$ encodes how elements of $\mathfrak{L}$ act as derivations on $A$, which is the essential ingredient for constructing a de Rham-type exterior algebra in this context.

\begin{definition}[Chevalley--Eilenberg algebra {\cite{rinehartDifferentialFormsGeneral1963}}]\label{def:CE-algebra}
Let $(\mathfrak{L}, \rho)$ be a Lie-Rinehart algebra over a commutative algebra $A$. The \emph{Chevalley--Eilenberg algebra} of $\mathfrak{L}$ is the graded $A$-algebra
\[
	\CE^\bullet(\mathfrak{L}) := 
    \left(\left(\Lambda_A^\bullet \mathfrak{L}\right)^\ast\right)^{\rev}~,
\]
where $\Lambda_A^\bullet \mathfrak{L} := S_A^\bullet(\mathfrak{L}[-1])$ denotes the graded-symmetric algebra over $A$ generated by the degree-shifted module $\mathfrak{L}[-1]$, the superscript ${*}$ stands for the $A$-linear dual $\Hom_A(\Lambda_A^\bullet \mathfrak{L}, A)$, and $\rev$ indicates reversal of the grading\footnote{The choice to reverse the grading is motivated by the convention that, when $\mathfrak{L}=\g$ is an ordinary Lie algebra (i.e., concentrated in degree $0$), one usually wants to consider both its exterior algebra $\Lambda\g$ and its Chevalley-Eilenberg complex $\CE^\bullet(\g)$ as non-negatively graded algebras, while taking the dual of a positively graded vector space yields a negatively graded one.}.
The algebra $\CE^\bullet(\mathfrak{L})$ is canonically endowed with the wedge product and a differential $\dCE \colon \CE^k(\mathfrak{L}) \to \CE^{k+1}(\mathfrak{L})$ defined by
\begin{align*}
	\dCE(\alpha)(X_1 \wedge \dotsb \wedge X_{k+1}) =\ & \sum_{i=1}^{k+1} (-1)^{i+1} \rho(X_i)\big(\alpha(X_1 \wedge \dotsb \widehat{X_i} \dotsb \wedge X_{k+1})\big) \\
	& + \sum_{1 \leq i < j \leq k+1} (-1)^{i+j} \alpha\big([X_i, X_j] \wedge X_1 \wedge \dotsb \widehat{X_i} \dotsb \widehat{X_j} \dotsb \wedge X_{k+1}\big)~,
\end{align*}
where $\alpha \in \CE^k(\mathfrak{L})$ is a $k$-form on $\mathfrak{L}$ and $X_1, \dotsc, X_{k+1} \in \mathfrak{L}$. The notation $\widehat{X_i}$ indicates omission of the $i$-th term from the wedge product.
\end{definition}

Differential forms on a smooth manifold provide a prototypical example.
\begin{example}[Classical de Rham algebra]
	When $\mathfrak L=\X( A)$, the corresponding Chevalley-Eilenberg algebra is called \emph{de Rham algebra} and it is denoted by $\Omega( A)$.
	If $A = C^\infty(M)$ for a smooth manifold $M$, then $\mathfrak{L} = \X(M)$ and $\big(\Omega(A), \wedge, \delta\big)$ is the classical de Rham algebra of smooth differential forms.
	A direct calculation, analogous to the smooth manifold case, shows that the CE-differential is, in fact, $A$-linear.
\end{example}

\begin{remark}
	A different approach to defining the de Rham complex over general commutative algebras, using universal properties, is described in \cite[\S 17]{Nestruev2020}. 
    \noindent For the constructions developed in this paper, the Lie--Rinehart–style formulation proves to be the more convenient choice.

\end{remark}

If we combine the above-defined $\dCE$ with the contraction operator $\iota:\mathfrak{L}\otimes \CE^\bullet(\mathfrak{L})\to \CE^\bullet(\mathfrak{L})$ defined by 
$$
	(\iota_X\alpha)(X_1,...,X_{k-1}):=\alpha(X,X_1,...,X_{k-1}) ~,
$$ 
and the Lie derivative 
$$
	\Lie_X\alpha=[\iota_X,\dCE]\alpha=\iota_X\dCE+\dCE\iota_X\alpha ~,
$$ 
the triple $(\dCE,\iota,\Lie)$ satisfies the standard rules of Cartan calculus
\begin{align*}
	\dCE^2                            &= 0,            &&&
	\iota_X \iota_Y + \iota_Y \iota_X  &= 0 ~,         \\
\iota_X \dCE   +	\dCE \iota_X      &= \Lie_X ~,    &&&
	\dCE \Lie_X - \Lie_X \dCE          &= 0 ~,    \\
	\Lie_X \iota_Y - \iota_Y \Lie_X    &= \iota_{[X,Y]}~,  &&&
	\Lie_X \Lie_Y - \Lie_Y \Lie_X      &= \Lie_{[X,Y]}~. 
\end{align*}
We will see in the next section how to reinterpret these operations as an action of a Gerstenhaber algebra on a BV-module.

\subsubsection{Gerstenhaber algebras}

Gerstenhaber algebras were introduced to formalise the algebraic structures arising in Hochschild cohomology and deformation theory, where a graded commutative product interacts compatibly with a Lie bracket of degree $-1$~\cite{Gerstenhaber1963}.

\begin{definition}[Gerstenhaber algebra]\label{def:Gerstenhaber-algebra}
	A \emph{Gerstenhaber algebra} is a triple $(G, \wedge, \{\cdot,\cdot\})$ where:
	\begin{itemize}
		\item $G$ is a $\mathbb{Z}$-graded vector space;
		\item $\wedge : G \otimes G \to G$ is a graded commutative, associative product;
		\item $\{\blank,\blank\} :G[1]\otimes G[1]\to G[1]$ is a graded Lie bracket,
	\end{itemize}
	such that, for any $a \in G$, the operation $\{a, -\}$ is a derivation of degree $|a|-1$ with respect to the product $\wedge$, that is,
	\begin{align}\label{eq:Gerstenhaber-derivationProperty}
		\{a, b \wedge c\} = \{a, b\} \wedge c + (-1)^{(|a|-1)|b|} b \wedge \{a, c\}
		&& \forall~ a, b, c \in G.
	\end{align}
Morphisms of Gerstenhaber algebras are maps which preserve both the associative and the Lie algebraic structure.
\end{definition}

\begin{remark} The bracket $\{\blank,\blank\}$ could also be reinterpreted as a degree $-1$ bracket on the desuspended $\Z$-graded vector space $G$, i.e. \autoref{def:Gerstenhaber-algebra} coicides with the definition of $1$-Poisson algebra of \cite{Cattaneo2006b}.
\end{remark}

Our primary example of a Gerstenhaber algebra is the exterior algebra of a Lie-Rinehart algebra.
\begin{lemma}\label{lem:LRyieldsGH}
Let $(A, \mathfrak{L})$ be a Lie--Rinehart algebra. Then the exterior algebra
\[
\Lambda_A \mathfrak{L} := S_A\left(\mathfrak{L}[-1]\right)
\]
acquires the structure of a Gerstenhaber algebra when endowed with the (associative) wedge product and the \emph{Schouten--Nijenhuis bracket}.
The Schouten--Nijenhuis bracket is the unique graded Lie bracket satisfying the following properties:
\begin{itemize}
	\item It restricts to the given Lie bracket on $\Lambda_A^1 \mathfrak{L} = \mathfrak{L}$;
	\item It vanishes on $A = \Lambda_A^0 \mathfrak{L}$;
	\item It satisfies the graded Leibniz rule with respect to the wedge product.
\end{itemize}
Explicitly, for homogeneous elements $x_1, \dots, x_m, y_1, \dots, y_n \in \mathfrak{L}$ and $f \in A$, it is given by
\begin{align*}
	[x_1 \wedge \cdots \wedge x_m,\; y_1 \wedge \cdots \wedge y_n] 
	&= \sum_{i,j} (-1)^{i+j} [x_i, y_j] \wedge x_1 \wedge \cdots \widehat{x_i} \cdots \wedge x_m \wedge y_1 \wedge \cdots \widehat{y_j} \cdots \wedge y_n, \\
	[f,\; x_1 \wedge \cdots \wedge x_m] 
	&= -\iota_{\dCE f}(x_1 \wedge \cdots \wedge x_m),
\end{align*}
where $\iota_{\dCE f}$ denotes contraction with the Chevalley--Eilenberg differential of $f$, and $\widehat{\,\blank\,}$ indicates omission.

This assignment forms a functor, i.e. Lie--Rinehart morphisms induce corresponding Gerstenhaber morphisms.

\end{lemma}

\begin{example}[Multivector fields]
	We have seen in \autoref{ex:XMisLR} that $A=C^\infty(M)$, $\mathfrak{L}=\mathfrak{X}(M)$ forms a Lie--Rinehart algebra. Hence, the above Lemma induces the classical Gerstenhaber algebra structure on $\X^\bullet(M):=\wedge_{C^\infty(M)} \X(M)\cong \bigoplus_k \Gamma(\wedge^k TM)$.
\end{example}

\begin{example}[Lie algebras]
	Since any Lie algebra $\g$ is a Lie--Rinehart algebra over the base field $\mathbb{k}$ with trivial anchor $\rho$, then the free graded commutative algebra $\Lambda\g$ admits a natural Gerstenhaber algebra structure, with the wedge product and the bracket given by the canonical extension of the Lie bracket (the Schouten bracket).
\end{example}

\subsubsection{Batalin–Vilkovisky modules}\label{sec:BVmodules}
In this section, we will recall Batalin–Vilkovisky (BV) modules as a special type of module over a Gerstenhaber algebra. BV-modules provide a very convenient way to extend classical Cartan calculus to a broad range of (commutative or non-commutative) algebraic settings \cite{Tsygan2004, Kowalzig2015}.

\begin{definition}[Gerstenhaber module]\label{def:GerstenhaberModule}
Let $(G,\wedge,\{\blank,\blank\})$ be a Gerstenhaber algebra and $V$ a graded vector space. A \emph{Gerstenhaber $G$-module structure} on $V$ consists of a Gerstenhaber algebra structure $(\star, \llbracket\blank,\blank\rrbracket)$ on the graded vector space $G \oplus V^{\rev}$ such that:
\begin{itemize}
	\item $G$ is a Gerstenhaber subalgebra of $(G \oplus V^{\rev}, \star, \llbracket\blank,\blank\rrbracket)$;
	\item Both brackets and associative multiplication vanish when restricted to $V^{\rev}$.
\end{itemize}
\end{definition}

\begin{remark}\label{rk:GerstenhaberModuleActions}
	We will denote the Gerstenhaber algebra structure on $G \oplus V^{\rev}$ by $(\star,\llbracket\blank,\blank\rrbracket)$, and the induced operations on $V$ by $\iota$ and $\Lie$.
	\\
	We will also denote the degree of an element $x \in G$ by $|x|$, and the degree of an element $\alpha \in V$ by $|\alpha|$.
	\\
	For a homogeneous element $x \in G$, we will denote by $\iota_x$ and $\Lie_x$ the corresponding operators on $V$.
The above definition encodes the standard notion of a module over a Gerstenhaber algebra as a graded vector space $V$ equipped with:
\begin{itemize}
	\item an associative action $\iota$ of $(G,\wedge)$ on $V$;
	\item a Lie action $\Lie$ of $(G,\{\blank,\blank\})$ on $V$ of degree $-1$.
\end{itemize}

These actions are induced by restricting the Gerstenhaber structure $(\star,\llbracket\blank,\blank\rrbracket)$ on $G \oplus V^{\rev}$ to $G \otimes V^{\rev} \subset (G \oplus V^{\rev})^{\otimes 2}$. We use the notations $\iota$ and $\Lie$ for these operations:
\begin{itemize}
	\item The associative action $\iota$ is defined by
	\[
	\morphism{\iota}{G \otimes V^{\rev}}{V^{\rev}}{(x,\alpha)}{\iota_x \alpha := x \star \alpha},
	\]
	and satisfies
	\[
	\iota_{x \wedge y} \alpha = \iota_x \iota_y \alpha \quad \text{and} \quad \iota_x \alpha = (-1)^{|x||\alpha|} \alpha \star x.
	\]

	\item The Lie action $\Lie$ is defined by
	\[
	\morphism{\Lie}{G \otimes V^{\rev}}{V^{\rev}[-1]}{(x,\alpha)}{\Lie_x \alpha := \llbracket x, \alpha \rrbracket},
	\]
	and satisfies
	\[
	\Lie_{\{x,y\}} \alpha = \Lie_x \Lie_y \alpha + (-1)^{(|x|-1)(|y|-1)} \Lie_y \Lie_x \alpha,
	\]
	as well as the graded antisymmetry identity
	\[
	\Lie_x \alpha = -(-1)^{|x||\alpha|} \llbracket \alpha, x \rrbracket.
	\]
\end{itemize}
\noindent
The compatibility between $\star$ and $\llbracket\blank,\blank\rrbracket$ in $G \oplus V^{\rev}$ (cf.~\autoref{eq:Gerstenhaber-derivationProperty}) implies the following \emph{mixed Leibniz rules} (see \cite[Def.~2.1]{Kowalzig2015}):
\begin{align}
	\iota_{\{x,y\}} \alpha &= \Lie_x \iota_y \alpha - (-1)^{(|x|-1)|y|} \iota_y \Lie_x \alpha, \\
	\Lie_{x \wedge y} \alpha &= \iota_x \Lie_y \alpha + (-1)^{|x||y|} \iota_y \Lie_x \alpha. \label{eq:mixedleib}
\end{align}
\noindent
The use of the reverse grading on $V$ in \autoref{def:GerstenhaberModule} ensures that for any homogeneous $x \in G^k$, the operators $\iota_x$ and $\Lie_x$ have degrees $-k$ and $1-k$, respectively, when viewed as endomorphisms of $V$.
\end{remark}

\begin{remark} \label{rk:metamodule}
	The above definition fits into a broader meta-concept: given a type of algebraic structure $\Gamma$ (e.g., associative, Poisson, Gerstenhaber, or Lie algebra) and a $\Gamma$-algebra $A$, the data of a module structure over $A$ can often be encoded as a $\mathbb{k}$-vector space $M$ together with a $\Gamma$-structure on $A \oplus M$, which:
	\begin{itemize}
		\item the inclusion $A \hookrightarrow A \oplus M$ is a $\Gamma$-algebra morphism,
		\item squares to zero, i.e., the structure is trivial when restricting everything to $M$.
	\end{itemize}
	This viewpoint, formalised via operads~\cite{KrizMay1995}, recovers the usual notions of modules for commutative and Lie algebras, and provides a natural framework for more intricate structures such as Gerstenhaber and Poisson modules. For more complex algebraic structures, multiple inequivalent notions of module may arise (cf. also \autoref{rk:othermods} below), but this construction typically yields the most structured version.
\end{remark}

The notations introduced in \autoref{rk:GerstenhaberModuleActions} for the actions of a Gerstenhaber algebra $G$ on a graded vector space $V$ are deliberately chosen to evoke the classical operations of Cartan calculus. 
To recover the full structure of Cartan calculus, one can equip $V$ with a differential as follows:

\begin{definition}[Batalin--Vilkovisky module]\label{def:BVmodule}
Let $G$ be a Gerstenhaber algebra, and let $(V, \d)$ be a cochain complex. 
A \emph{Batalin--Vilkovisky (BV) $G$-module} structure on $V$ consists of a Gerstenhaber $G$-module structure together with a differential $d$ of degree $+1$ such that, for all $x \in G$,
\begin{align}\label{eq:Cartanmagic}
\Lie_x = [\iota_x, d] = \iota_x \circ d - (-1)^{|x|} d \circ \iota_x.
\end{align}
\end{definition}
\begin{remark}\label{rk:retrieve-multivector-Cartan-Calculus}
	In this setting, the operators $\d$, $\iota_x$, and $\Lie_x$ satisfy six graded commutation relations of classical Cartan calculus (extended to multivector fields):
	\begin{align*}
		\Lie_x&=[\iota_x,d]\\
		\Lie_{[x,y]}&=[\Lie_x,\Lie_y]\\
		\iota_{[x,y]}&=[\Lie_x,\iota_y]=(-1)^{|y|}[\iota_x,\Lie_y]\\
		\iota_{xy}&=\iota_x\iota_y\\
		\Lie_{xy}&=\iota_x\Lie_y + (-1)^{|y|}\Lie_x \iota_y
	\end{align*}
	A pair $(G, V)$ consisting of a Gerstenhaber algebra and a BV module is often called a \emph{Tsygan–Tamarkin–Nest noncommutative differential calculus}~\cite[Def.~4.4]{Tsygan2004}.
\end{remark}

Uniting the Gerstenhaber algebra from \autoref{lem:LRyieldsGH} and the Chevalley-Eilenberg algebra from \autoref{def:CE-algebra}, we obtain a natural class of BV-modules.	

\begin{example}[BV modules from Lie--Rinehart algebras]
\label{ex:BVmodulesfromLieRinehart}
	Let $\mathfrak{L}$ be a Lie-Rinehart algebra over a commutative algebra $A$. This structure naturally gives rise to Gerstenhaber and Batalin--Vilkovisky (BV) algebraic objects, generalising classical Cartan calculus~\cite{Huebschmann1998}.

\begin{itemize}
    \item As discussed in \autoref{lem:LRyieldsGH}, the exterior algebra $\Lambda_A^\bullet \mathfrak{L}$, equipped with the wedge product $\wedge$ and the Schouten--Nijenhuis bracket, forms a Gerstenhaber algebra.

    \item The Chevalley--Eilenberg algebra $\CE(\mathfrak{L}) := \left((\Lambda_A^\bullet \mathfrak{L})^*\right)^{\rev}$, introduced in \autoref{def:CE-algebra}, carries a canonical structure of Gerstenhaber module over $\Lambda_A^\bullet \mathfrak{L}$.\footnote{The two grading reversals from \autoref{def:CE-algebra} and \autoref{def:GerstenhaberModule} cancel out, i.e. we have a Gerstenhaber algebra structure on $\Lambda\mathfrak{L}\oplus(\Lambda\mathfrak{L})^*$. However, if one omits both reversals, the formulas and gradings one obtains become inconsistent with the literature.} The module operations on an element $\alpha \in \CE(\mathfrak{L})$ are defined as follows:
    \begin{itemize}
        \item For $a \in A$, we set $\iota_a \alpha := a \cdot \alpha$. For $X \in \mathfrak{L}$, the contraction $\iota_X \alpha$ is the usual insertion into the leftmost slot. For elements of higher wedge degree, the contraction is extended recursively by the rule\footnote{Some sources, e.g.~\cite{rogers2012linfty,CalliesFregierRogersZambon16, Ryvkin2018, blacker2023reduction}, use the convention $\tilde\iota_{X_1 \wedge \dots \wedge X_k} \alpha := \iota_{X_k} \dots \iota_{X_1} \alpha$, which differs from ours by a sign $(-1)^{\frac{k(k-1)}{2}}$.}:
        \[
        \iota_{X_1 \wedge \dots \wedge X_k} \alpha := \iota_{X_1} \dots \iota_{X_k} \alpha.
        \]

        \item The Lie derivative is defined via the Cartan formula:
        \[
        \Lie_X := \iota_X \circ \dCE - (-1)^{|X|} \dCE \circ \iota_X.
        \]
    \end{itemize}

    \item Equipped with the Chevalley--Eilenberg differential $\dCE$ and the operators $\iota$ and $\Lie$, the complex $\CE(\mathfrak{L})$ acquires the structure of a BV module over $\Lambda_A^\bullet \mathfrak{L}$.
\end{itemize}
\end{example}

\begin{example}[Classical Cartan calculus as a BV module]\label{ex:BVmodulesfromCartan}
	Let $M$ be a smooth manifold. The Lie--Rinehart algebra $\mathfrak{X}(M)$ of vector fields over the commutative algebra $C^\infty(M)$ induces, via the construction of \autoref{ex:BVmodulesfromLieRinehart}, a canonical Batalin--Vilkovisky module structure on the de Rham complex $\Omega^\bullet(M)$.

	Here, the Gerstenhaber algebra is given by the exterior algebra \( G = \Lambda^\bullet \mathfrak{X}(M) \) of (positively graded) multivector fields, endowed with the wedge product and the Schouten--Nijenhuis bracket.

	The corresponding BV module is the complex \( \CE(\mathfrak{X}(M)) := \left((\wedge^\bullet \mathfrak{X}(M))^*\right)^{\rev} \), which identifies naturally with the (positively graded) space of differential forms \( \Omega^\bullet(M) \). The structure maps are:
	\begin{itemize}
		\item the de Rham differential $\d = \dCE$,
		\item the contraction $\iota_X$ by multivector fields $X \in G$, defined by insertion into the leftmost slot,
		\item the Lie derivative $\Lie_X$ via Cartan's formula \( \Lie_X = [\iota_X, \d] \).
	\end{itemize}
	This recovers the classical Cartan calculus as a Batalin-Vilkovisky module structure on the de Rham complex, governed by the Gerstenhaber algebra of multivector fields.
\end{example}

\begin{remark}[Variants of Gerstenhaber modules]\label{rk:othermods}
	The notion of a Gerstenhaber module is not uniquely fixed across the literature. Besides differing sign conventions, a key point of variation lies in whether the so-called mixed Leibniz rule of \autoref{eq:mixedleib} is imposed explicitly (cf.~\cite{Kowalzig2014,DeSole2011,Tsygan2004}). 
	For BV modules, however, the mixed Leibniz rule follows automatically from Cartan-type identities and thus does not need to be postulated separately.
\end{remark}

\subsection{{$L_\infty$-algebras} from BV-modules}
\label{sec:linftyfrombv}
The construction of the $L_\infty$-algebra in \autoref{def:linftygeom} can be made general. It can be carried out for any BV-module equipped with a fixed $k+1$-cocycle \( \omega \), i.e., for which \( \d \omega = 0 \) with respect to the differential \( \d \), since being a BV-module provides the structure of a complete Cartan calculus, which is all we need. In particular, the following definition generalises the notion of Hamiltonian pairs and the corresponding observable algebra.

\begin{definition}[Hamiltonian pairs]\label{def:algebraicham0}
	Let \( (G, \wedge, \{\blank, \blank\}) \) be a Gerstenhaber algebra, and \( (V, \d) \) a BV-module over \( G \).
	Assume \( \omega \in V^{k+1} \) is a $k+1$-cocycle, i.e., \( \d \omega = 0 \) for $k\geq 1$.
	The space of \emph{Hamiltonian pairs} is defined as the graded vector space:
	\[
		\Ham^0(V, \omega) := \left\{ (\alpha, X) \in V^{k-1} \oplus G^1 \mid \iota_X \omega = -\d\alpha \right\},
	\]
	which can also be viewed as the fibered product \( V^{k-1} \times_{V^k} G^1 \) over the following pullback in the category of vector spaces:
	\begin{displaymath}
		\begin{tikzcd}
			\Ham^0(V, \omega) \arrow[r] \arrow[d] \ar[dr,"\lrcorner",very near start,phantom]& G^1 \ar[d,"\iota_{\blank}\omega"] \\
			V^{k-1} \arrow[r, "-\d"] & V^{k}
		\end{tikzcd}
	\end{displaymath}
\end{definition}

We can now extend the definition of $L_\infty$-algebras, as in  \autoref{def:linftygeom}, to this algebraic context. This leads to the following construction of the (generalised) $L_\infty$-algebra of observables.

\begin{definition}[Generalized Observables $L_\infty$-algebra]\label{def:AlgebraicRogersAlgebra}
	Let \( (G, \wedge, \{ \blank, \blank \}) \) be a Gerstenhaber algebra, and \( (V, \d) \) a BV-module over \( G \).
	Fix a $k+1$-cocycle \( \omega \in V^{k+1} \) with \( \d\omega = 0 \), for $k\geq 1$.
	We define the \emph{graded space of algebraic observables} as
	$$
		\Ham(V, \omega) := \bigoplus_{i\leq 0} \Ham^i(V, \omega),
	$$
	where
	$$
		\Ham^0(V, \omega) := \left\{ (\alpha, X) \in V^{k-1} \oplus G^1 \mid \iota_X \omega = -\d\alpha \right\}
	$$
	is the space of algebraic Hamiltonian pairs (\autoref{def:algebraicham0}), and for each $ i \leq -1$,
	$$
		\Ham^i(V, \omega) := V^{k-1 + i}~.
	$$

	The graded vector space $\Ham(V, \omega)$ carries an $L_\infty$-algebra structure with degree $(2-j)$ multibrackets $l_j\colon \Ham^{\otimes j}(V, \omega) \to \Ham(V, \omega)[2-j]$, defined as:
	\begin{itemize}
		\item \( l_1 \colon \Ham^i(V, \omega) \to \Ham^{i+1}(V, \omega) \):
		      \[
			      l_1(\alpha) = \begin{cases}
				      \d\alpha      & \text{if } i < -1, \\
				      (\d\alpha, 0) & \text{if } i = -1, \\
				      0             & \text{if } i = 0;
			      \end{cases}
		      \]
		\item \( l_2 \colon \Ham^0(V, \omega) \otimes \Ham^0(V, \omega) \to \Ham^{0}(V,\omega) \):
		      \[
			      l_2\big( (\alpha, X), (\beta, Y) \big) = \big( \iota_X \iota_Y \omega, \{X, Y\} \big);
		      \]
		\item For each $j \geq 3$~:
		      \[
			      l_j\big( (\alpha_1, X_1), \dots, (\alpha_j, X_j) \big) = -\iota_{X_1} \dots \iota_{X_j} \omega;
		      \]
		\item All multibrackets of arity greater than or equal to $2$ vanish when evaluated on inputs containing at least one element of nonzero degree.
	\end{itemize}

	We refer to \( \Ham(V, \omega) \) endowed with these brackets as the \emph{$L_\infty$-algebra of observables}.
\end{definition}
	The verification that the above definition yields an honest $L_\infty$-algebra can be directly adapted from~\cite[Thm. 5.2]{rogers2012linfty} noticing that their proof is purely algebraic in nature and relies only on the Cartan calculus axioms.

\begin{example}[Observables $L_\infty$-algebra associated with a Lie--Rinehart algebra]\label{rem:AlgebraicRogersAlgebra}
	According to \autoref{ex:BVmodulesfromCartan}, the construction of the $L_\infty$-algebra in \autoref{def:AlgebraicRogersAlgebra} can be carried out in the context of any Lie–Rinehart algebra \( (\mathfrak{L}, A) \), equipped with a Chevalley–Eilenberg $(k+1)$-cocycle \( \omega \) for the differential \( \dCE \).
	Specifically, one obtains $L_\infty$-brackets on the complex
	\begin{displaymath}
		\begin{tikzcd}[]
			A \arrow[r] &
			\mathfrak{L}^* \arrow[r] &
			\cdots \arrow[r] &
			(\Lambda^{k-2} \mathfrak{L})^* \arrow[r] &
			(\Lambda^{k-1} \mathfrak{L})^* \times_{(\Lambda^k \mathfrak{L})^*} \mathfrak{L}
		\end{tikzcd}
	\end{displaymath}
	where the fibre product is defined by the same relation as in \autoref{def:algebraicham0}, and $(\Lambda^k \mathfrak{L})^*$ corresponds to the Chevalley--Eilenberg algebra of \autoref{def:CE-algebra} in degree $k$ (as such, the corresponding BV-module is concentrated in non-negative degrees).
	We note that in \cite{richterLieInfinityAlgebrasLie2013}, the author proposes a distinct $L_\infty$-algebra associated with a Lie--Rinehart algebra which is not directly related to the one we are constructing here.
\end{example}

\begin{example}\label{ex:clasfits}
	In particular, for $A=C^\infty(M)$ and $\mathfrak{L}=\mathfrak{X}(M)$, \autoref{def:linftygeom} becomes an example of \autoref{def:AlgebraicRogersAlgebra}.
\end{example}

\begin{remark} We refer to \cite[\S 4]{Delgado2018b} for a slightly different interpretation of the  $L_\infty$-algebra of multisymplectic observables using Gerstenhaber algebras, including the construction of a bigger $L_\infty$-algebra in the setting of \autoref{ex:clasfits}.
\end{remark}

\section{Constraint Cartan Calculus}\label{sec:constraint-Cartan-Calculus}
In this section, we review the general theory of constraint algebras and modules, and introduce a notion of \emph{constraint Cartan calculus} adapted to BV-modules within this framework.

\subsection{Reminder on Constraint Triples}
We begin by recalling the algebraic theory of constraint triples, originally introduced in~\cite{dippelespositowaldmann}. Our exposition follows the more detailed treatment in the PhD thesis~\cite{Dippell2023}, which serves as our primary reference. For simplicity, we work over a fixed base field $\mathbb{k}$ of characteristic zero, although most of the constructions extend to more general settings.

\subsubsection{Constraint Vector Spaces}
We start by introducing the basic notion of a constraint vector space, which underlies the formalism of constraint triples. This provides the linear algebraic foundation upon which more elaborate constraint structures are built.
\begin{definition}[Contraint vector spaces] An \emph{(embedded, $\mathbb Z$-graded) constraint vector space} $V$ is given by a triple $V=(V_\Total,V_\Wobs,V_\Null)$, where $V_\Total$ is a $\mathbb Z$-graded vector space and $V_\Null\subset V_\Wobs\subset V_\Total$ are subspaces. A \emph{morphism} $f:V\to V'$ of constraint vector spaces is a linear map $f:V_\Total\to V'_\Total$ such that $f(V_\Wobs)\subset V'_\Wobs$ and $f(V_\Null)\subset V'_\Null$.
\end{definition}
\noindent We think of such a constraint triple to be a "pre-reduction" description of the subquotient $\frac{V_\Wobs}{V_\Null}$. 

\begin{remark}[Non-embedded constraint objects]\label{rk:embed}
    In this article, we will only treat embedded objects, i.e., cases where $V_\Wobs$ is a subspace of $V_\Total$. 
    In general, the constraint formalism can also handle cases where the inclusion of $V_\Wobs$ into $V_\Total$ is replaced by a general map $V_\Wobs\to V_\Total$.  
\end{remark}

\begin{remark}[Reduction functor]\label{rk:red-functor}
	The very purpose of introducing categories of constraint objects is to assume the existence of a corresponding \emph{reduction functor}, typically denoted by $\red$. 
	For example, in the category of constraint vector spaces, the reduction functor $\red$ assigns to each object $V$ the quotient vector space $\red(V) := V_\Wobs/V_\Null$, and to each morphism $f: V \to V'$ the induced map $\red(f): V_\Wobs/V_\Null \to V'_\Wobs/V'_\Null$.
\end{remark}

\begin{definition}\label{def:homotimesR} Let $V$ and $V'$ be constraint vector spaces. 
	\begin{enumerate}
		\item The direct sum $V\oplus V'$ is defined as the triple $(V_\Total\oplus V'_\Total,V_\Wobs\oplus V'_\Wobs,V_\Null\oplus V'_\Null)$. It is the product and coproduct of the category of constraint vector spaces.
		\item The internal hom $\iHom(V,V')$ is defined as the triple $$\Triple{\iHom(V,V')}=
		\colvec{
		\Hom(V_\Total,V'_\Total)\\
		\{f\in\Hom(V_\Total,V'_\Total)~|~f(V_\Wobs)\subset V'_\Wobs,f(V_\Null)\subset V'_\Null \} \\
		\{f\in\Hom(V_\Total,V'_\Total)~|~f(V_\Wobs)\subset V'_\Null \} 		
		}$$ 
		\item The (weak) tensor product of $V$ and $V'$ is defined as	
	$$
	\Triple{(V\otimes V')}=\colvec{
		V_\Total\otimes V'_\Total\\
		 V_\Wobs\otimes V'_\Wobs\\
		 V_\Null\otimes V'_\Wobs + V_\Wobs\otimes V'_\Null}
	$$
	\end{enumerate}
\end{definition}
We note that the usual homomorphisms from $V$ to $V'$ in the category of constraint vector spaces are exactly elements of $\iHom(V, V')_\Wobs$.  The definition of the tensor product can be justified by the following Lemma (which is proven in the ungraded case in \cite{Dippell2023}, but whose proof works analogously in the graded world):

\begin{lemma}{\cite[Propositions 1.2.22, 1.2.23]{Dippell2023}}
\begin{enumerate}
\item The tensor product $\otimes$ turns constraint vector spaces into a monoidal category with monoidal unit $\underline{\mathbb{k}}:=(\mathbb{k},\mathbb{k},0)$.
\item The reduction functor $V=\triple{V}\mapsto \frac{V_\Wobs}{V_\Null}$ from constraint vector spaces to vector spaces is monoidal.
\item For a fixed constraint vector space $V$ the functor $\cdot \otimes V$ is left adjoint to $\iHom(V,\cdot)$ in the category of constraint vector spaces.
\end{enumerate}

\end{lemma}

\subsubsection{Constraint algebras and their modules}
%
With the monoidal structure on the category of constraint vector spaces in place, we can naturally extend the classical notion of associative algebras to this setting. This leads to the definition of \emph{constraint algebras}, along with their associated modules, which form the basic algebraic structures of interest in our framework.
\begin{definition}[Constraint associative algebra]
	A \emph{constraint associative algebra} is a monoid $(A,m:A\otimes A\to A,u: \underline{\mathbb{k}}\to A)$ in the monoidal category of constraint vector spaces.
\end{definition}

\begin{remark}[Constraint associative algebra as a triple]\label{rk:constraintalgebras-as-triples}
	Unpacking the definitions, a constraint associative algebra is equivalently given by a triple $A = \triple{A}$, where $A_\Wobs \subset A_\Total$ is a subalgebra, and $A_\Null \subset A_\Wobs$ is a (two-sided) ideal. In particular, this ensures that the reduced space $\red(A) := A_\Wobs / A_\Null$ inherits a natural associative algebra structure, making the reduction functor $\red$ (see \autoref{rk:red-functor}) compatible with the algebraic operations.
\end{remark}

Graded commutative algebras can be obtained from associative ones by requiring that the multiplication (on $A_\Total$) is graded-commutative. Similarly, we can define constraint Lie algebras.

\begin{definition}[Constraint Lie algebra]
	A \emph{constraint (graded) Lie algebra} is a triple $\g=\triple{\mathfrak g}$, where $\g_\Total$ is a (graded) Lie algebra, $\g_\Wobs\subset \g_\Total$ is a sub-Lie-algebra and $\g_\Null\subset \g_\Wobs$ a Lie ideal.
\end{definition}

\begin{example}\label{ex:cinfycalg} Let $M$ be a smooth manifold and $S\subset M$ a closed set. We write $I_S\subset C^\infty(M)$ for the ideal of smooth functions vanishing on $S$. Then $A=(C^\infty(M),C^\infty(M),I_S)$ is a commutative constraint algebra. 
\end{example}

\begin{example}\label{ex:homsalgcstr}
The ($\mathbb{k}$-linear) endomorphisms $\iHom(A,A)$ of $A$ form an associative constraint algebra when equipped with composition and a constraint Lie algebra when equipped with the commutator.
\end{example}

\begin{example}\label{ex:geometricconstrainttriplealgebrazied}
	The algebraic formulation of symplectic reduction recalled in the introduction, see in particular \autoref{eq:algredexchange}, can be naturally interpreted as a constraint algebra considering the following triple:
	\[
	\Triple{{A}} := 
	\colvec{
	C^\infty(M_{\Total}) = C^\infty(M) \\
		 \{ f \in \mathcal{A}_{\Total} \mid \forall g \in G: \vartheta_g^*f - f \in {A}_\Null \} \\
		 I_{M_\Wobs} = I_{\mu^{-1}(0)}
	}
	\]
	One easily verifies that ${A}_\Wobs$ is a subalgebra: for any $f_1, f_2 \in {A}_\Wobs$ and $g \in G$, we compute
	\begin{align*}
		\vartheta_g^*(f_1 f_2) - f_1 f_2 
		&= \vartheta_g^*(f_1)\vartheta_g^*(f_2) - f_1 f_2 \\
		&= \vartheta_g^*(f_1)(\vartheta_g^*(f_2) - f_2) + (\vartheta_g^*(f_1) - f_1)f_2 \in {A}_\Null.
	\end{align*}
	The triple $({A}_{\Total}, {A}_\Wobs, {A}_\Null)$ thus provides a prototypical example of an associative \emph{constraint algebra}.
\end{example}

Using the technique described in \autoref{rk:metamodule}, we can directly obtain the notion of modules over these types of algebras.

\begin{definition}[Constraint algebra module] Let $A$ be a constraint algebra (graded-commutative or Lie) and $\mathcal{M}$ a constraint vector space. An \emph{$A$-module structure} on $\mathcal{M}$ is given by a constraint algebra structure of the same type on $A\oplus \mathcal{M}$, such that $A$ is a subalgebra and the multiplication is trivial on $\mathcal{M}$. 
\end{definition}

\begin{remark}[Embedded setting]\label{rk:simpler-def-on-Total}
	If we unpack this definition, the module structure is specified by a constraint morphism $A \otimes \mathcal{M} \to \mathcal{M}$, such that the induced map on the $T$-components endows $M_\Total$ with the structure of an $A_\Total$-module of the given type. This strategy naturally extends to the definition of other types of constraint algebraic structures by prescribing the relevant operation as a multilinear constraint morphism and imposing the defining axioms a priori only on the $ T$-component. This is possible because all of our objects are \emph{embedded} (cf.~\autoref{rk:embed}), which implies that the forgetful functor mapping the constraint object to its $T$-component is faithful, i.e. morphisms are entirely determined by what they do on $T$-components.
	
\end{remark}

As a refinement of the $\mathbb{k}$-linear endomorphism Lie algebra from \autoref{ex:homsalgcstr}, we can also consider derivations.

\begin{example}[The constraint space of derivations]\label{ex:cderismod} Let $A$ be a commutative constraint algebra concentrated in degree zero. 
	We define the derivations of $A$ as the triple
	$$
	 \mathfrak X(A)=\colvec{
	 	\mathfrak{X}(A_\Total)\\
	 	\iHom(A,A)_\Wobs\cap 	\mathfrak{X}(A_\Total)\\
	 	\iHom(A,A)_\Null\cap 	\mathfrak{X}(A_\Total)
	 }
	$$
	I.e. we are taking those elements of $\iHom$ which satisfy the derivation properties. $\mathfrak{X}(A)$ is both an $A$-module and a Lie algebra.
\end{example}

Now that we have the notion of modules over constraint algebras, we can try to mimic \autoref{def:homotimesR} in the $A$-linear (rather than just ${\mathbb{k}}$-linear case). For homomorphisms, everything works as expected; however, for tensor products, we have to be a little careful:

\begin{definition} Let $A$ be a graded-commutative constraint algebra and $\mathcal{M},\mathcal{M}'$ constraint $A$-modules.
\begin{enumerate}
	\item We define the $A$-linear internal Hom $A$-module to be given by:
	$$
	\iHom_A(\mathcal{M},\mathcal{M}')=\Triple{\iHom_A(\mathcal{M},\mathcal{M}')}=
	\colvec{
		\iHom(\mathcal{M},\mathcal{M}')_\Total\cap \Hom_A(\mathcal{M}_\Total,\mathcal{M}'_\Total)\\
		\iHom(\mathcal{M},\mathcal{M}')_\Wobs\cap \Hom_A(\mathcal{M}_\Total,\mathcal{M}'_\Total)\\
		\iHom(\mathcal{M},\mathcal{M}')_\Null\cap \Hom_A(\mathcal{M}_\Total,\mathcal{M}'_\Total)
	}
	$$
	\item In particular, the dual $A$-module $\mathcal{M}^*$ is defined as $\iHom_A(\mathcal{M},A)$.
	\item The tensor product $\mathcal{M}\otimes_A\mathcal{M}'$ is defined as:
	$$
\mathcal{M}\otimes_A\mathcal{M}'=\Triple{(\mathcal{M}\otimes_A\mathcal{M}')}=
\colvec{
	\mathcal{M}_\Total\otimes_{A_\Total} \mathcal{M}'_\Total\\
    j(\mathcal{M}_\Wobs\otimes_{A_\Wobs} \mathcal{M}'_\Wobs)\\
    j(\mathcal{M}_\Null\otimes_{A_\Wobs} \mathcal{M}'_\Wobs+ \mathcal{M}_\Wobs\otimes_{A_\Wobs} \mathcal{M}'_\Null)
},
$$
where $j$ is the natural map $\mathcal{M}_\Wobs\otimes_{A_\Wobs} \mathcal{M}'_\Wobs\to \mathcal{M}_\Total\otimes_{A_\Total} \mathcal{M}'_\Total$.
\end{enumerate}
\end{definition}

\begin{remark} The presence of $j$ in the definition of the tensor product is there to assure that $(\mathcal{M}\otimes_A\mathcal{M}')_\Wobs$ is a subspace of $(\mathcal{M}\otimes_A\mathcal{M}')_\Total$, i.e. that we don't leave the embedded setting (cf. \autoref{rk:embed}). A non-embedded version of the tensor product is discussed in \cite{Dippell2023}.
\end{remark}

\begin{example}\label{ex:constextalg}
	Let $A$ be a commutative constraint algebra concentrated in degree 0 and $\mathcal{M}$ an $A$-module also concentrated in degree zero. By considering the graded symmetrisation of $\bigotimes_A^k\mathcal{M}[1]$ (in all components) and summing over all $k$, we obtain the exterior algebra $\Lambda_A\mathcal{M}$. In components, it reads:
	$$
		\Lambda_A\mathcal{M}=
		\colvec{
			\Lambda_{A_\Total}\mathcal{M}_\Total\\
			j(\Lambda_{A_\Wobs}\mathcal{M}_\Wobs)\\
			j(\mathcal{M}_\Null\wedge {\Lambda_{A_\Wobs}}\mathcal{M}_\Wobs)
		}
	$$ 
where this time $j$ is the natural map from $\Lambda_{A_\Wobs}\mathcal{M}_\Wobs$ to $\Lambda_{A_\Total}\mathcal{M}_\Total$.
\end{example}

\subsubsection{Strong Constraint algebras}
%
In practice, it is relatively common for the $A_\Null$ component to not only be an ideal in $A_\Wobs$ but also in $A_\Total$. Geometrically, in the setting of \autoref{eq:algredexchange}, this corresponds to the group action being defined on all of the smooth manifold $M$, not only on the zero level set of the corresponding momentum map. We hence define

\begin{definition}[Strong constraint algebra]
	A constraint algebra $A=\triple{A}$ is called strong if $A_\Null$ is an ideal in $A_\Total$.
\end{definition}

In order to understand strong constraint algebras as monoids in a monoidal category, one can make the following definition.

\begin{definition}[Strong tensor product]
	Let $V,V'$ be constraint vector spaces, then their strong tensor product is defined as:
	$$
	V\boxtimes V'=\Triple{(V\boxtimes V')}=
	\colvec{
		V_\Total\otimes V'_\Total\\
		V_\Wobs\otimes V'_\Wobs + 
		V_\Null\otimes V'_\Total + V_\Total\otimes V'_\Null\\
		V_\Null\otimes V'_\Total + V_\Total\otimes V'_\Null}
	$$
\end{definition}
Strong constraint algebras are monoids in the monoidal category of constraint vector spaces with $\boxtimes$ as the monoidal structure. We will provide a slightly different characterisation, using the fact that
there is a natural constraint homomorphism $j: V\otimes V'\to V\boxtimes V'$. This homomorphism allows us to formulate the following
\begin{lemma} A constraint algebra $A$ is strong if its multiplication $m:A\otimes A\to A$ factors through $A\boxtimes A$, i.e. there exists $\overline m: A\boxtimes A\to A$ such that $m=\overline m\circ j$.
\end{lemma}

One can also define a strong $A$-linear tensor product, but we omit its definition, as it will not be needed in the sequel. 

Many algebraic constructions admit strong versions, for instance:

\begin{definition}[Strong constraint module] 
	Let $A$ be a strong constraint graded-commutative algebra. An $A$-module $\mathcal{M}$ is called strong if the module multiplication $A\otimes \mathcal{M}\to \mathcal{M}$ factors through $A\boxtimes \mathcal{M}$.
\end{definition}

\begin{remark} Dualization does not preserve strongness (or non-strongness) of structures. For instance, for finite-dimensional constraint vector spaces (which are automatically strong as $\underline{k}$-modules), we have $(V\otimes V')^*=V^*\boxtimes (V')^*$ and $(V\boxtimes V')^*=V^*\otimes (V')^*$.
\end{remark}

\begin{example}\label{ex:constextalgdual}
	Let $\mathcal{M}$ be a strong $A$-module for a strong constraint algebra $A$. Let $\Lambda_A\mathcal{M}$ be the exterior algebra from \autoref{ex:constextalg}. In general, $\Lambda_A\mathcal{M}$ is not strong as an algebra. However, the constraint algebra structure on its dual $(\Lambda_A\mathcal{M})^*$ is strong. We can look at its components as follows:
	
$$
(\Lambda_A^k\mathcal{M})^*=
\colvec{
(\Lambda^k_{A_\Total}\mathcal{M}_\Total)^*\\
\{
f\in (\Lambda_{A_\Total}\mathcal{M}_\Total)^* | ~ f(v_1,...,v_k)\in A_\Wobs,~ f(w,v_2,...,v_k)\in A_\Null~\forall v_1,...,v_k\in \mathcal{M}_\Wobs, \forall w\in \mathcal{M}_\Null
\}\\
\{
f\in (\Lambda_{A_\Total}\mathcal{M}_\Total)^* | ~ f(v_1,...,v_k)\in A_\Null~\forall v_1,...,v_k\in \mathcal{M}_\Wobs
\}
}
$$	
To see that it is strong, one has to verify that for $f\in (\Lambda_A^k\mathcal{M})^*_\Null$ and $g\in (\Lambda_A^l\mathcal{M})^*_\Total$, the product $f\wedge g$ lies in $(\Lambda_A^{k+l}\mathcal{M})^*_\Null$. 
\end{example}

We close the subsection by defining (strong) constraint Lie--Rinehart algebras.

\begin{definition}[Constraint Lie--Rinehart algebra] Let $A$ be a constraint commutative algebra. A constraint Lie Rinehart algebra over $A$ is an $A$-module $\mathfrak{L}$ together with a constraint Lie algebra structure and an $A$-linear constraint morphism $\rho:\mathfrak{L}\to \mathfrak X(A)$, such that $\mathfrak{L}_\Total$ is a Lie--Rinehart algebra over $A_\Total$. We call the pair $A,\mathfrak{L}$ strong when $A$ is a strong constraint algebra and $\mathfrak{L}$ is a strong constraint module.
\end{definition}  

Of course, $\mathfrak{X}(A)$ itself is a constraint Lie--Rinehart algebra for any constraint algebra $A$ (cf. \autoref{ex:cderismod}). When $A$ is strong, it is even a strong Lie--Rinehart algebra.

We are not requiring the Lie bracket of $\mathfrak{L}$ to be strong, because in all cases relevant to us it will not be. This can be illustrated by the following example.
\begin{example}\label{ex:cinftyctgt}Let $M$ be a smooth manifold and $S$ a closed subset. We denote by $I_S$ the ideal of all functions vanishing on $S$. As in \autoref{ex:cinfycalg}, we can now consider the strong constraint algebra 
	$$
	A=(C^\infty(M),C^\infty(M),I_S)~.
	$$
Its derivations are:
$$
\mathfrak X(A)=\colvec{
\mathfrak X(M)\\
\{X\in\mathfrak X(M) ~|~X(I_S)\subset I_S\}\\
\{X\in\mathfrak X(M) ~|~X(C^\infty(M))\subset I_S\}
}
$$
We think of $\mathfrak{X}(A)_\Wobs$ to be the space of vector fields (algebraically) tangent to $S$\footnote{We use the adjective algebraically here, because it can be a distinct notion from being geometrically tangent, cf. \autoref{ex:algtang} later in the text.} and  $\mathfrak{X}(A)_\Null$ are the vector fields vanishing on $S$. The Lie bracket of $ \mathfrak{X}(A)$ is not strong, even in the simplest cases: If $M=\mathbb R$ and $S=\{0\}$, then $x\partial_x\in  \mathfrak{X}(A)_\Null$ but $[\partial_x,x\partial_x]=\partial_x\not\in  \mathfrak{X}(A)_\Null$ contradicting strongness of the bracket.
\end{example}

\subsection{Constraint BV-modules}\label{sec:constraintBVmodules}
We have seen in \autoref{sec:BVmodules} how the notion of BV-modules can be understood as a Cartan calculus. In this section, we will adapt this notion to the setting of constraint algebras and modules.
\begin{definition}[Constraint Gerstenhaber algebra {(\cite[Def. 2.25]{dippelkern})}]\label{def:constraintgerst}
	Let $G=\triple{G}$ be a constraint algebra. We say that $G$ is a \emph{constraint Gerstenhaber algebra} if it is equipped with a constraint bilinear map
	\begin{displaymath}
		\{ \blank, \blank \} \colon G[1] \otimes G[1] \to G[1]
	\end{displaymath}
	such that $G_\Total$ is a Gerstenhaber algebra as per \autoref{def:Gerstenhaber-algebra}.
	A constraint Gerstenhaber algebra 
	${G} $ is called 
	\emph{strong}, if the underlying graded constraint algebra is strong. 
\end{definition}

 \begin{remark}\label{rk:constraintgerst-as-triple}
	Since we are in the embedded setting (see \autoref{rk:embed}), a constraint Gerstenhaber $G$ as per \autoref{def:constraintgerst} can be seen as a triple of Gerstenhaber subalgebras $G_\Null \subseteq G_\Wobs \subseteq (G_\Total, \wedge,\lbrace\blank,\blank\rbrace)$ such that $G_\Null$ is a Gerstenhaber ideal in $G_\Wobs$. 
	I.e., the following four equations hold:
	\begin{align*}
		\{ G_\Wobs, G_\Wobs \} &\subseteq G_\Wobs~,&& G_\Wobs \wedge G_\Wobs \subseteq G_\Wobs~, \\
		\{ G_\Null, G_\Wobs \} &\subseteq G_\Null~,&& G_\Wobs \wedge G_\Null \subseteq G_\Null ~.
	\end{align*}
	For a strong Gerstenhaber algebra, we would also need to check $G_\Null \wedge G_\Total \subseteq G_\Null$ as the only additional condition.
\end{remark}

Since every constraint Gerstenhaber algebra is, in particular, a graded constraint algebra and a graded constraint Lie algebra, we know that ${G}_{\red}$ carries the structures of a graded algebra and a graded Lie algebra.
Moreover, it is easy to see that the graded Leibniz rule on ${G}$
induces a graded Leibniz rule on ${G}_{\red}$, showing that ${G}_{\red}$ forms itself a Gerstenhaber algebra.

\begin{definition}[Constraint Gerstenhaber module]\label{def:constraintgerstmodule}
	Let $G=\triple{G}$ be a constraint Gerstenhaber algebra. 
	A \emph{constraint $G$-module} is a constraint vector space $V=\triple{V}$ endowed with a constraint Gerstenhaber algebra structure on the constraint vector space $G \oplus V^\rev$ such that the multiplication and the bracket are trivial when restricted to $V^\rev$.
\end{definition}
\autoref{def:constraintgerstmodule} implies that $G_\Total \oplus V^\rev_\Total$ is a Gerstenhaber algebra and, according to the spirit of \autoref{rk:constraintgerst-as-triple}, a constraint $ G$-module can be seen as a sequence of Gerstenhaber submodules $V_\Null \subseteq V_\Wobs \subseteq V_\Total$.

Furthermore, one can introduce the notion of constraint BV-modules employing the strategy of \autoref{rk:simpler-def-on-Total}.
\begin{definition}[Constraint BV-module]\label{def:constraintbv}
	Let $G=\triple{G}$ be a constraint Gerstenhaber algebra. A \emph{constraint BV-module} is a constraint Gerstenhaber module $V=\triple{V}$ endowed with a constraint degree $1$ homomorphism
	$$
		\d \colon V \to V[1]~,
	$$
	such that $(V_\Total, \d_\Total)$ is a BV-module over $G_\Total$ as per \autoref{def:BVmodule}.
\end{definition}

\begin{remark}
	Spelling out \autoref{def:constraintbv} in the same way of \autoref{rk:GerstenhaberModuleActions} one can reframe a constraint BV-module structure on $V$ as the existence of three constraint morphisms $\iota: G\otimes V^{\rev} \to V^{\rev}$, $\d: V \to V[1]$ and $\Lie: G\otimes V^{\rev} \to V^{\rev}[-1]$ satisfying six commutation rules analogous to \autoref{rk:retrieve-multivector-Cartan-Calculus}.
	Therefore, one can understand a constraint BV-module as a \emph{constraint Cartan calculus}.
\end{remark}

Similarly to the non-constraint case (see \autoref{ex:BVmodulesfromLieRinehart}), one can construct a constraint BV-module out of any constraint Lie--Rinehart algebra.
\begin{lemma}\label{lem:const-bv-from-const-LR}
	Let $\mathfrak L$ be a constraint Lie--Rinehart algebra over a constraint algebra $A$. 
	Then the constraint vector space $\Lambda_A \mathfrak L$ is a constraint Gerstenhaber algebra and $\CE(\mathfrak L) = ((\Lambda_A \mathfrak{L})^*)^\rev$ is a constraint BV-module over it.
\end{lemma}
\begin{proof}
	The explicit expression of $\Lambda_A \mathfrak L$ and $(\Lambda_A \mathfrak{L})^*$ as contraint triples has been already given in \autoref{ex:constextalg} and \autoref{ex:constextalgdual}, respectively. By unpacking these definitions, one can also verify that the contraction $\iota$ is a constraint morphism. What remains to be shown is that $\dCE$ naturally defined on the total components according to \autoref{ex:BVmodulesfromLieRinehart} descends to a proper constraint homomorphism on the triple. Then the same will hold true automatically for $\Lie$ by Cartan's magic formula. 
	\begin{itemize}\item Let $\alpha\in \CE^k(\mathfrak L)_\Null$. We have to show that $\dCE\alpha$ is also in the Null component, i.e. that $$\dCE(\alpha)(X_1,...,X_{k+1})\in A_\Null$$ when for all $i$ we have $X_i\in \mathfrak{L}_\Wobs$. This follows from the definition of $\dCE$, recalling that $[\cdot,\cdot]$ preserves $\mathfrak{L}_\Wobs$, $\alpha$ maps tuples from $\mathfrak{L}_\Wobs$ to $A_\Null$ and $\rho(X_i)(A_\Null)\subset A_\Null$.
	\item For $\alpha\in \CE^k(\mathfrak L)_\Wobs$ we have to show two things, firstly that
	$\dCE(\alpha)(X_1,...,X_{k+1})\in A_\Wobs$  when for all $i$ we have $X_i\in \mathfrak{L}_\Wobs$. and secondly that $\dCE(\alpha)(X_1,...,X_{k+1})\in A_\Null$ if additionally $X_1\in\mathfrak{L}_\Null$. Both follow from the same type of argument by using that the Lie bracket and $\rho$ are constraint morphisms, with the only additional subtlety of a case distinction in the last calculation, depending on whether $\rho$ is applied to $X_1$ or $X_i$ ($i\neq 1$) in the second sum defining $\dCE$.
 	\end{itemize}
\end{proof}

\begin{remark}
The Cartan calculus developed in~\cite[\S3.4]{dippelkern} can be seen as a special case of the construction presented here. We will explore this in \autoref{ss:geomcartan}.
	%
	%
\end{remark}

\begin{remark}
	When a BV-module $V$ arises as the Chevalley--Eilenberg algebra of a Lie--Rinehart algebra, it inherits the structure of a strong algebra with respect to the wedge product~$\wedge$.
	In contrast, the noncommutative differential calculus developed by Tsygan, Tamarkin, and Nest (see \autoref{rk:retrieve-multivector-Cartan-Calculus}) does not, in general, require a wedge product structure.
\end{remark}

\section{Multisymplectic constraint reduction}\label{sec:multisymplecticconstraintreduction}
%
In this section, we will extend the construction of an $L_\infty$-algebra associated with a BV-module, as proposed in \autoref{def:AlgebraicRogersAlgebra}, to the context of constraint algebras developed in \autoref{sec:constraintBVmodules}. 
Our goal is to associate a BV-module and a corresponding $L_\infty$-reduction with any given pair of associative algebra and ideal. We will then apply this framework to the case of the Hamiltonian $L_\infty$-algebra of a multisymplectic manifold in the presence of symmetries.

\subsection{Constraint $L_\infty$-algebras and their reduction}
%
Building on the chosen embedded constraint setting (see \autoref{rk:embed}), one can further define constraint $L_\infty$-algebras in close analogy with \autoref{def:constraintbv}, as follows.

\begin{definition}[Constraint $L_\infty$-algebras]\label{def:constraintlinfty}
	Let $L=\triple{L}$ be a constraint graded vector space.
	We say that $L$ is a \emph{constraint $L_\infty$-algebra} if it is equipped with a collection of constraint multilinear maps
	\begin{displaymath}
		l_j \colon L^{\otimes j} \to L[2-j]
	\end{displaymath}
	such that $(L_\Total, (l_j)_\Total)$ is an $L_\infty$-algebra (see \cite{Lada1993} for the original formulation and \cite{Ryvkin2016a} for a modern overview relevant to our setting).
\end{definition}
 \begin{remark}[Constraint $L_\infty$-algebra as a triple]\label{rk:constLinfty-as-triple}
	Since we work in the embedded setting (compare with \autoref{rk:constraintgerst-as-triple}), a constraint $L_\infty$-algebra $L$ as per \autoref{def:constraintlinfty} is given by a sequence of $L_\infty$-subalgebras $L_\Null \subseteq L_\Wobs \subseteq (L_\Total,(l_j)_\Total)$ such that $L_\Null$ is in particular an $L_\infty$-ideal  in $L_\Wobs$ in the following sense:
	\begin{displaymath}
		l_j\big(L_\Null, L_\Wobs,\dots, L_\Wobs\big)
		\subseteq L_\Null~, \quad \forall j \geq 1~.
	\end{displaymath}
	Therefore, is assured that the $L_\infty$-algebra structure on $L_\Total$ induces a well-defined $L_\infty$-algebra structure on the reduced space $\red(L) := L_\Wobs / L_\Null$ (compare with \autoref{rk:constraintalgebras-as-triples}).
\end{remark}

We now extend the construction of \autoref{sec:linftyfrombv} to the constraint setting, i.e. we define a constraint $L_\infty$-algebra of Hamiltonian pairs out of any constraint BV-module with a suitably fixed cocycle.
\begin{definition}[Contraint vector space of Hamiltonian pairs]\label{def:constraint-algebraicham0}
	Let $G=\triple{G}$ be a constraint Gerstenhaber algebra and $V=\triple{V}$ a constraint BV-module over $G$.
	Let $\omega\in V_\Wobs^{k+1}$ be a closed $(k+1)$-form, i.e. a $(k+1)$-cocycle with respect to $\d$. 
	The space of \emph{Hamiltonian pairs} is defined as the constraint graded vector space:
	\[
		\Ham^0(V, \omega) := 
		\colvec{
			\big\{ (\alpha, X) \in V^{k-1}_\Total \oplus G_\Total^1 \mid \iota_X \omega = -\d\alpha \big\}
			\\
			\Ham^0(V, \omega) \cap \big(V_\Wobs \oplus G_\Wobs\big)
			\\
			\Ham^0(V, \omega) \cap \big( V_\Null \oplus G_\Null \big)
		}
		~.
	\]
\end{definition}
\begin{remark}[Hamiltonian pairs as a fibered product]
	The space of Hamiltonian pairs can be interpreted as the fibered product \( V^{k-1} \times_{V^k} G^1 \) over the same diagram as in \autoref{def:algebraicham0}, now considered within the category of constraint vector spaces. 
	This perspective yields, in particular, three pullback diagrams, one for each component of the constraint vector space described in \autoref{def:constraint-algebraicham0}.
\end{remark}
	The constraint vector space of Hamiltonian pairs can be enlarged to give rise to a constraint $L_\infty$-algebra, as follows.
\begin{lemma}\label{lem:clinftyconstr}
	Let \( (G, \wedge, \{ \blank, \blank \}) \) be a constraint Gerstenhaber algebra, and \( (V, \d) \) a constraint BV-module over \( G \).
	Fix a $k+1$-cocycle \( \omega \in V^{k+1}_\Wobs \) with \( \d\omega = 0 \).
	The \emph{constraint graded vector space}
	\[
		\Ham(V, \omega) := \bigoplus_{i \leq 0} \Ham^i(V, \omega),
	\]
	given in degree $0$ by the constraint vector space of Hamiltonian pairs \( \Ham^0(V, \omega) \) as per \autoref{def:constraint-algebraicham0}, and by
	\[
		\Ham^i(V, \omega) := V^{k-1 + i} \qquad ~, \forall i \leq -1
	\]
	forms a constraint $L_\infty$-algebra when endowed with the family of constraint multilinear maps \( l_j \colon \Ham^{\otimes j}(V, \omega) \to \Ham(V, \omega)[2-j] \) defined as in \autoref{def:AlgebraicRogersAlgebra}.
\end{lemma}
\begin{proof}
	By its very definition, $\Ham(V, \omega)_\Total$ is an $L_\infty$-algebra since it is given specializing \autoref{def:AlgebraicRogersAlgebra} to the BV-module $V_\Total$ over $G_\Total$. Thus, we only have to check that the brackets are constraint morphisms. For this, we consider $l_1$ and the remaining cases separately:
	\begin{itemize}
		\item Since $V$ is a constraint BV-module, $l_1:\Ham^{i}(V, \omega)\to \Ham^{i+1}(V, \omega)$ is a constraint morphism for $i<-1$. For $i=1$, the same holds true because the additional component $0$ in the expression $l_1(\alpha)=(d\alpha,0)$ lies in $V_\Null$.
		\item The higher brackets $l_j$, $j\geq 3$ are constraint morphisms, since the projections $\Ham^0(V,\omega)\to G^1$ and the contractions $G^1\times V\to V$ are. For $l_2$, we need to additionally use the fact that the Lie bracket $G^1\otimes G^1 \to G^1$ is a constraint morphism.
	\end{itemize}

\end{proof}
	We refer to \( \Ham(V, \omega) \) endowed with these brackets as the \emph{$L_\infty$-algebra of observables}.

\subsection{Constraint $L_\infty$-algebra induced by symmetries}
\label{sec:constraintlinftyred-derivations}

In this section, we are going to carry out our main construction to obtain a constraint $L_\infty$-algebra of observables associated to an appropriate constraint algebra $A$, a collection of symmetries $F$ and a cocycle $\omega$. \\

%
Throughout, we fix a strong constraint commutative algebra \( A \) satisfying \( A_\Total = A_\Wobs \). (In \autoref{sec:constraintlinftyred-multisymp}, we will specialize to the case where \( A_\Wobs \) is the algebra of smooth functions on a manifold, and \( A_\Null \) is the ideal of functions vanishing on a closed subset, as in \autoref{ex:cinfycalg}.)
The corresponding strong constraint Lie--Rinehart algebra of derivations is given by
$$
\mathfrak{X}(A) :=
\colvec{
	\mathfrak{X}(A_\Total) \\
	\{X \in \mathfrak{X}(A_\Total) \mid X(A_\Null) \subset A_\Null\} \\
	\{X \in \mathfrak{X}(A_\Total) \mid X(A_\Wobs) \subset A_\Null\}
}.
$$
Following \autoref{ex:cinftyctgt}, we interpret $\mathfrak{X}(A)_\Wobs$ as the space of vector fields (algebraically) preserving the subset defined by the ideal $A\Null$, and $\mathfrak{X}(A)_\Null$ as the subspace of vector fields vanishing on that subset.
Let now $\mathcal{F}$ be a Lie--Rinehart algebra over $A_\Wobs$ satisfying the inclusions $\mathfrak{X}(A)_\Null \subset \mathcal{F} \subset \mathfrak{X}(A)_\Wobs$. We regard $\mathcal{F}$ as specifying a choice of symmetries of the system (e.g. it may arise from the infinitesimal action of a Lie algebra on a smooth manifold) preserving the closed subset corresponding to $A_\Null$.
This setting allows us to define an algebraic analogue of the constraint triple introduced in \autoref{ex:geometricconstrainttriplealgebrazied}.

\begin{lemma}The triple 
	$$
	A'=\Triple{A'}=\colvec{
	A_T\\
	\{f\in A_\Wobs~|~ \mathcal F(f)\subset A_0\}\\
	A_\Null
	}
	$$
	defines a strong constraint algebra.
\end{lemma}
\begin{proof} 
	According to \autoref{rk:constraintalgebras-as-triples}, it suffices to verify that $A'_\Wobs$ is a subalgebra of $A_\Total$ and that $A'_\Null$ is an ideal in $A'_\Wobs$.
	The only nontrivial condition is that $A'_\Wobs$ is closed under multiplication, i.e. $A'_\Wobs\cdot A'_\Wobs\subset A'_\Wobs$. To verify it, let $f,g\in A'_\Wobs$ and $X\in \mathcal F$. We have to show $X(fg)\in A'_0$. For this, we calculate:
	$$
	X(fg)=fX(g)+g X(f)\in A'_\Wobs\cdot A'_\Null + A'_\Wobs\cdot A'_\Null\subset A'_\Null
	~.
	$$
\end{proof}

By construction, $A'_\Wobs$ consists of those elements of $A_\Wobs$ that are mapped into $A_\Null = A'_\Null$ by the action of elements of $\mathcal{F}$. However, there may exist additional derivations in $\mathfrak{X}(A_\Total)$ that satisfy the same property. We denote the space of all such derivations by $\overline{\mathcal{F}}$:
\begin{equation}\label{eq:barF}
	\overline{\mathcal{F}} := \{ X \in \mathfrak{X}(A_\Total) \mid X(A'_\Wobs) \subset A'_\Null \}
	\supseteq \mathcal{F}~.
\end{equation}
This maximal space of derivations moving $A'_\Wobs$ into $A'_\Null$ naturally arises when studying the constraint space of derivations of the constraint algebra $A'$, as formalised in the following lemma.

\begin{lemma} \label{lem:Aprimeconstr}
	The constraint Lie--Rinehart algebra of derivations $\mathfrak{X}(A')$ is given by
	$$
	\mathfrak{X}(A')=\colvec{
	\mathfrak X(A_\Total)\\
	\mathcal{N}(\overline{\mathcal{F}},\mathfrak{X}(A)_\Wobs)\\
	\overline{\mathcal{F}}
	}
	$$
	where $\mathcal{N}(\overline{\mathcal{F}},\mathfrak{X}(A)_\Wobs) := \lbrace X \in \X(A)_\Wobs ~\vert~ [\overline{\mathcal{F}},X]\subseteq \overline{\mathcal{F}}\rbrace$ denotes the Lie algebra normalizer of $\overline{\mathcal{F}}$ in $\mathfrak{X}(A)_\Wobs$.
\end{lemma}

\begin{proof}
	The first and last components follow directly from \autoref{ex:cderismod} and the definition of $\overline{\mathcal{F}}$ in \autoref{eq:barF}, respectively. 
	For the middle component, we need to identify the normalizer $\mathcal{N}(\overline{\mathcal{F}}, \mathfrak{X}(A)_\Wobs)$ with the space of derivations $X$ satisfying both $X(A'_\Null) \subset A'_\Null$ and $X(A'_\Wobs) \subset A'_\Wobs$. 
	The first condition is automatically fulfilled for all elements of $\mathfrak{X}(A)_\Wobs$, as $A'_\Null = A_\Null$ by construction. 
	To verify the second condition, we proceed with the following computation:
\begin{align*}
	&X(A'_\Wobs)\subset A'_\Wobs&\\
	\iff & Y(X(A'_\Wobs))\subset A_\Null &\forall Y\in\overline{\mathcal{F}}\\
	\iff & (X\circ Y - [X,Y])(A'_\Wobs)\subset A_\Null &\forall Y\in\overline{\mathcal{F}}
\end{align*}	
Since $Y\in \overline{\mathcal{F}}$ we have $Y(A'_\Wobs)\subset A_\Null$, hence $X\circ Y(A'_\Wobs)\subset	A_\Null$, so we can omit this term and obtain:
\begin{align*}
	\iff &  [X,Y](A'_\Wobs)\subset A_\Null &\forall Y\in\overline{\mathcal{F}}\\
	\iff &  [X,Y]\in\overline{\mathcal{F}} &\forall Y\in\overline{\mathcal{F}}\\
\end{align*}	
But the last statement is precisely the defining equality for the normaliser in question.
\end{proof}

We now return to our original choice of $\mathcal{F}$ and proceed to construct a corresponding Lie-Rinehart algebra.
\begin{proposition}\label{prop:Yconstr} 
	The following triple defines a strong Lie--Rinehart algebra over $A'$:
	$$
	\mathfrak{Y}=\colvec{
	\mathfrak{X}(A)_\Total\\
	\mathcal{N}({\mathcal{F}},\mathfrak{X}(A)_\Wobs)\\
	\mathcal{F}
	}~.
	$$
\end{proposition}

\begin{proof}The Lie bracket preserves the constraint structure due to the properties of the normaliser. So let us consider the $A'$-module structure:
\begin{itemize}
\item We start by showing that $A'_\Wobs\mathfrak Y_\Wobs\subset Y_\Wobs$. Let $f\in A'_\Wobs$ and $X\in Y_\Wobs$. First of all $$fX\in A'_\Wobs\mathfrak{X}(A)_\Wobs\subset A_\Wobs\mathfrak{X}(A)_\Wobs\subset \mathfrak{X}(A)_\Wobs.$$
Now let $Y\in \mathcal F$, we want to show that $[Y, fX]\in \mathcal F$. We have 
$$
[Y,fX]=Y(f)X+f[X,Y]\in Y(A'_\Wobs)X+ A'_\Wobs [X,Y]\subset A'_\Null X + A'_\Wobs \mathcal F\subset \mathfrak X(A)_\Null + \mathcal{F}\subset \mathcal{F}
$$
This means that $fX\in \mathfrak{Y}_\Wobs$.
\item The inclusion $A'_\Null\mathfrak{Y}_\Total\subset \mathfrak{Y}_\Null$ is a consequence of $A_\Null\mathfrak X(A)_\Total\subset\mathfrak{X}(A)_\Null\subset \mathcal{F}$. 
\item The inclusion $A'_\Total\mathfrak{Y}_\Null\subset \mathfrak{Y}_\Null$ is a consequence of the fact that $\mathcal F$ is a Lie--Rinehart algebra over $A_\Total$.
\end{itemize}
What remains to show is that $\mathfrak{Y}$ acts by derivations of $A'$, i.e. that the identity induces a constraint morphism $\rho:\mathfrak{Y}\to \mathfrak{X}(A')$. Since the total components agree, and on the null components we have $\mathcal{F}\subset\overline{\mathcal{F}}$, we only have to show that
$$
    \mathcal{N}({\mathcal{F}},\mathfrak{X}(A)_\Wobs)\subset \mathcal{N}(\overline{\mathcal{F}},\mathfrak{X}(A)_\Wobs)~.
$$

Let $X\in \mathfrak{X}(A)_\Wobs$ be given, such that $[X,\mathcal{F}]\subset\mathcal{F}$. We have to show $[X,\overline{\mathcal{F}}]\subset\overline{\mathcal{F}}$. We first observe that $X(A'_\Wobs)\subset A'_\Wobs$. To verify this, one can pick $f\in  A'_\Wobs$ and verify that $X(f)\in  A'_\Wobs$, which means that for all  $Y\in\mathcal{F}$, we have $Y(X(f))\in A_\Null$. To do this, we calculate:
$$
Y\circ X(f)=X\circ Y(f)+[X,Y](f)\in X(A_\Null) + \mathcal F(f)\subset A_\Null.
$$
Now let $Z\in \overline F$. We can calculate:
\begin{align*}
[X,Z](A'_\Wobs)&\subset X\circ Z (A'_\Wobs) -  Z\circ X (A'_\Wobs)\\
&\subset X(A_\Null) - Z (A'_\Wobs) \subset A_\Null - A_\Null\subset A_\Null.
\end{align*}
Hence, $[X,Z]\in\overline{\mathcal{F}}$, which concludes the proof.
\end{proof}

The Lie--Rinehart algebra $\mathfrak{Y}$ should be seen as a refined version of $\mathfrak{X}(A')$, which remembers more about the original choice of symmetries. Multiple different $\mathcal{F}$s can have the same $\overline{\mathcal F}$. 
In the geometric setting, as discussed in \autoref{sec:constraintlinftyred-multisymp}, this additional memory is typically related to vanishing orders and flat points. 

\begin{remark}\label{rem:symgenLR} 
	In geometric contexts, symmetries are often given by a collection of vector fields—i.e., derivations—which might not form a Lie-Rinehart algebra. In those case, to an arbitrary subset $F \subset \mathfrak{X}(A)_\Wobs$, one can associate the minimal Lie--Rinehart algebra $\mathcal{F}_F$ containing both $F$ and $\mathfrak{X}(A)_\Null$, and perform the constructions above using $\mathcal{F}_F$ in place of $\mathcal{F}$.
\end{remark}

We now interpret the constraint Lie--Rinehart algebra $(A', \mathfrak{Y})$ as a \emph{subquotient} of $(A, \mathfrak{X}(A))$: that is, it can be obtained by first restricting to the subalgebra and then passing to a quotient. Concretely, this corresponds to the following diagram which should be read from right to left.
\[
\Triple{\mathfrak{X}(A)}
\overset{\mathrm{quotient}}\longleftarrow
\colvec{
	\mathfrak{X}(A)_\Total \\
	\mathfrak{Y}_\Wobs \\
	\mathfrak{X}(A)_\Null
}
\overset{\mathrm{subspace}}\longrightarrow
\Triple{\mathfrak{Y}}.
\]
We now wish to carry out an analogous construction on the dual side: namely, to define a procedure for first restricting and then projecting differential forms on $A$ to obtain differential forms on $A'$. Formally, this takes the following shape:

\begin{proposition} \label{prop:Bprimeconst}Let $B$ be a BV-module for $\Lambda_A\mathfrak X(A)$. Then the following triple constitutes a BV-module over $\Lambda_{A'}\mathfrak{Y}$
$$
B'=\Triple{B'}=\colvec{
B_\Total\\
\{\alpha\in B_\Wobs~|~\iota_{\mathfrak{Y}_\Null}\alpha, \Lie_{\mathfrak{Y}_\Null}\alpha \subset B_\Null \}\\
B_\Null
}.
$$
  
\end{proposition}
\begin{proof} 
	We start by noting that there is no canonical constraint Lie--Rinehart morphism from $\mathfrak{Y}$ to $\mathfrak{X}(A)$, and also no Gerstenhaber morphism between their exterior algebras. However, since the total spaces of $\mathfrak{Y}$,  $\mathfrak{X}(A)$ and $B, B'$ are the same, we can make sense of all operations; we just have to check that they verify the constraint structure.
\begin{itemize}
	\item We first take care of the differential $\d$ being a constraint morphism when restricted to $B'$. Essentially, we only have to check that $\d B'_\Wobs\subset B'_\Wobs$. This is the standard computation for basic forms in differential geometry forming a subcomplex. Let $\alpha\in B'_\Wobs$ and $X\in \mathfrak{Y}_\Wobs$. We can verify:
	\begin{align*}
		\Lie_X\d\alpha&=\d\Lie_X\alpha\in \d B_\Null\subset B_\Null\\	
		\iota_X \d\alpha&=\Lie_X\alpha - \d\iota_X\alpha\in B_\Null - \d B_\Null\subset B_\Null
	\end{align*}
	\item Contractions with elements of $\mathfrak{Y}_\Wobs$ preserve $B'_\Wobs$ by construction. So do contractions with $A'_\Wobs$: Let $a\in A'_\Wobs, \alpha\in B'_\Null, X\in \mathfrak{Y}_\Wobs$, then
	\begin{align*}
		\iota_X(a\alpha)&=a\iota_X\alpha\in a\cdot B_\Null\subset A'_\Wobs B_\Null\subset A_\Wobs B_\Null\subset B_\Null\\
		\Lie_X(a\alpha)&=a\Lie_X\alpha + X(a)\alpha	
		\in  B_\Null + A_\Null B_\Total\subset B_\Null
	\end{align*}
	Since the exterior algebra is generated by the degrees zero and one, this already implies that $\iota$ restricts to a morphism $(\Lambda_{A'}\mathfrak{Y})_\Wobs\otimes B'_\Wobs\to B'_\Wobs$. 
	\item  Similarly, we have $A'_\Null B'_\Wobs\subset A_\Null B_\Wobs\subset B_\Null$. The condition  $\iota_{\mathfrak{Y}_\Null}B'_\Wobs\subset B'_\Null$ follows from the definition of $B'_\Wobs$. The space  $(\Lambda_{A'}\mathfrak{Y})_\Null$ is generated by products $X\wedge \xi$ with $X\in \mathfrak{Y}_\Null$ and $\xi\in (\Lambda_{A'}\mathfrak{Y})_\Wobs$. Together with the above item, this implies that $\iota$ also restricts to a map $(\Lambda_{A'}\mathfrak{Y})_\Null\otimes B'_\Wobs\to B'_\Null$.
	\item The last property we need of the contraction is that it restricts to a map $(\Lambda_{A'}\mathfrak{Y})_\Wobs\otimes B'_\Null\to B'_\Null$. Since $(\Lambda_{A'}\mathfrak{Y})_\Wobs$ is generated by $A'_\Wobs$ and $\mathfrak{Y}_\Wobs$, we can again only show the statement for those types of elements. We have:
	\begin{align*}
		A'_\Wobs B'_\Null&\subset A_\Wobs B_\Null\subset B_\Null\\
		\iota_{\mathfrak{Y}_\Wobs} B'_\Null&\subset \iota_{\mathfrak{X}(A)_\Wobs} B'_\Null \subset B'_\Null.
	\end{align*} 
	\item Finally, the constraint identities for $\Lie$ follow from those for $\iota$ and those for $d$, due to \autoref{eq:Cartanmagic}.
\end{itemize}	
\end{proof}

Summarising, the above construction can be reformulated in terms of the following data:
\begin{enumerate}
	\item An (ordinary) commutative algebra $\overline A$;
	\item An ideal $I \subset \overline A$;
	\item A collection of symmetries preserving $I$, i.e., a subset $F \subset \mathfrak{X}(\overline A, I) := \{ X \in \mathfrak{X}(\overline A) \mid X(I) \subset I \}$.
\end{enumerate}

From this data, we constructed the strong constraint algebra $A := (\overline A, \overline A, I)$, together with its Lie--Rinehart algebra of derivations $\mathfrak{X}(A)$. 
Then, focused on the Lie--Rinehart algebra $\mathcal{F} := \mathcal{F}_F$ generated by $F + \mathfrak{X}(A)_\Null$ (cf.~\autoref{rem:symgenLR}), and used it to define the constraint algebra $A'$ (cf.~\autoref{lem:Aprimeconstr}) and the associated constraint Lie--Rinehart algebra $\mathfrak{Y}$ (cf.~\autoref{prop:Yconstr}).
We then considered the constraint BV-module $B := \mathrm{CE}(\mathfrak{X}(A))$ over the exterior algebra $\Lambda_A \mathfrak{X}(A)$, and the induced module of \emph{basic forms} $B'$ over $\mathfrak{Y}$ (cf.~\autoref{prop:Bprimeconst}).

Moreover, if we are given a cocycle $\omega \in \CE^{k+1}(\mathfrak{X}(A))_\Wobs$ such that 
$\iota_{\mathcal{F}}\omega,~\Lie_{\mathcal{F}}\omega \subset \CE(\mathfrak{X}(A))_\Null$, 
then $\omega$ descends to a cocycle in $B'$, yielding a constraint $L_\infty$-algebra as described in \autoref{lem:clinftyconstr}.
Remarkably, the condition on the Lie derivative is not necessary for this construction to work; it is already sufficient to require that 
$\iota_{\mathcal{F}}\omega \subset \CE(\mathfrak{X}(A))_\Null$. 
It is enough to check this property on the generating set $F \subset \mathcal{F}$. That is, we can equivalently formulate the requirement as follows:
\begin{enumerate}[resume]
	\item A closed element $\omega \in \CE^{k+1}(\mathfrak{X}(\overline{A}))$ (i.e., $\dCE \omega = 0$) satisfying
	\[
	\omega\big(F, \mathfrak{X}(\overline{A}, I), \ldots, \mathfrak{X}(\overline{A}, I)\big) \subset I.
	\]
\end{enumerate}
Putting everything together, we arrive at the following general recipe:

\begin{theorem}[Reduced algebra of observables]\label{thm:mainredgenalgleotired} 
	Given the data 1. - 4. above, we obtain a constraint $L_\infty$-algebra on the space $\Ham(V, \omega)$ as defined in \autoref{lem:clinftyconstr}, in particular an $L_\infty$-algebra structure on $\frac{\Ham(B', \omega)_\Wobs}{\Ham(B', \omega)_\Null}$. 
	We will call this algebra the \emph{reduced algebra of observables}.
\end{theorem}
\begin{proof} 
	There are two things we still need to prove: On the one hand, that the $L_\infty$-brackets preserve the constraint structure (even though under the weaker condition 4. $\omega$ need not restrict to a cocycle on $B'$), and the fact that $\iota_F\omega\subset \CE(\mathfrak{X}(A))_\Null$. 

	Let us start with the latter one.
	$\mathcal F$ is the minimal Lie--Rinehart algebra containing $F$ and $\mathfrak{X}(A)_0$. 
	$\iota_F\omega\subset  \CE(\mathfrak{X}(A))_\Null$ by 4., and $\iota_{\mathfrak{X}(A)_0}\omega\subset  \CE(\mathfrak{X}(A))_\Null$ since $\omega\in \CE(\mathfrak{X}(A))_\Wobs$. 
	The minimal Lie--Rinehart algebra is constructed from its generators by a series of module multiplications and Lie brackets. For module multiplications, this is a consequence of $CE(\mathfrak{X}(A))$ being a strong constraint algebra. 
	Hence, we can just start with $X,Y\in \mathfrak{X}(A)_\Wobs$ such that $\iota_X\omega,\iota_Y\omega \in \CE(\mathfrak{X}(A))_\Null$ and show that the same is true for the bracket. This goes as follows:
	$$
	\iota_{[X,Y]}\omega=\Lie_X\iota_Y\omega - \iota_Y\Lie_X\omega =
	\Lie_X\iota_Y\omega - \iota_Y\d\iota_X\omega\in \Lie_X \CE(\mathfrak{X}(A))_\Null + \iota_Y \d~ \CE(\mathfrak{X}(A))_\Null\subset \CE(\mathfrak{X}(A))_\Null
	$$
	Here we used that $X,Y\in \mathfrak X(A)_\Null$ so that $\Lie_X$, $\iota_Y$ and $\d$ preserve $\CE(\mathfrak{X}(A))_\Null$.
	
	Let us now turn to why the $l_k$ are constraint morphisms. 
	$l_1$ and the $\mathfrak{Y}$-valued part of $l_2$ do not depend explicitly on $\omega$, hence they are automatically constraint morphisms. 
	Thus, we only need to see why for $k\geq 2$ the maps
	\begin{align*}
		\Ham^0(B',\omega)_\Total\times  ....\times \Ham^0(B',\omega)_\Total\to B'_\Total\\
		(\alpha_1,X_1), ..., (\alpha_k,X_k)\mapsto -\iota_{X_1}...\iota_{X_k}\omega
	\end{align*}
	induce constraint morphisms. 
	However, $-\iota_{X_1}...\iota_{X_k}\omega$ can be rewritten as $\iota_{X_1}...\iota_{X_{k-1}}d\alpha$, which now does not depend on $\omega$ and induces a constraint morphism because $B'$ is a BV-module. 
\end{proof}

\subsection{$L_\infty$-reduction of multisymplectic manifolds}\label{sec:constraintlinftyred-multisymp}
	In this section, we apply the general, purely algebraic reduction procedure developed above to the geometric setting. 
	Namely, let $M$ be a smooth manifold and $S \subset M$ a closed subset, we specialise the construction of \autoref{sec:constraintlinftyred-derivations} to the case $\overline{A} = C^\infty(M)$ and $I = I_S$, the vanishing ideal of $S$. 
	This corresponds to considering the constraint algebra $A = (C^\infty(M), C^\infty(M), I_S)$, as in \autoref{ex:cinfycalg}, together with its Lie--Rinehart algebra of derivations $\mathfrak{X}(A)$, introduced in \autoref{ex:cinftyctgt}:
	$$
	\mathfrak X(A)=\colvec{
		\mathfrak X(M)\\
		\{X\in\mathfrak X(M) ~|~X(I_S)\subset I_S\}\\
		\{X\in\mathfrak X(M) ~|~X(C^\infty(M))\subset I_S\}
	}~.
	$$
	We first remark that $\mathfrak{X}(A)_0$ consists precisely of the vector fields vanishing on the set $S$. 
	This is because $X(f)=\d f(X)$ and since the possible $\d f$'s span the whole cotangent space, the only way how $X(C^\infty(M))\subset I_S$ can be satisfied is if $X|_S=0$. 
	In fact, we can even write such a vector field as a \underline{finite} sum of functions vanishing on $S$ times vector fields. This is a more general phenomenon:
	\begin{lemma}\label{lem:vanonsubset}
		Let $E$ be a vector bundle over $M$ and $S$ a closed subset. Then any section $s: M\to E$ vanishing on $S$ can be written as a finite sum of functions vanishing on $S$ with sections of $E$.
		In particular, in the above setting, $\mathfrak{X}(A)_\Null=I_S\cdot \mathfrak X(M).$
	\end{lemma}
	\begin{proof}
		We prove the first statement. For trivializable vector bundles, the statement is clear: one picks a frame $e_1,...,e_n$ of $E$ and the coefficient in front of the $e_j$ must be zero on $S$. 
		\\
		For a general vector bundle $E$ there is always another vector bundle $F$ such that $E\oplus F$ is trivial (see e.g. \cite[\S 12.33]{Nestruev2020}), i.e. admits a frame $e_1,...,e_N$. 
		Seeing a section $s$ of $E$ as a section of $E\oplus F$ we can find functions $f_i$ vanishing on $S$ such that $s=\sum_{i=1}^N f_i e_i$. 
		Now, by construction $s=\sum_{i=1}^Nf_i\pi_E(e_i)$ as a section of $E$.
\end{proof}

The above Lemma has direct consequences for constraint differential forms:
\begin{corollary}\label{cor:bnbt}
	The constraint BV-module $B=\CE(\mathfrak{X}(A))$ satisfies $B_\Wobs=B_\Total$.
\end{corollary}
\begin{proof}
Let $\alpha\in B_\Total$, we want to show that it automatically lies in $B_\Wobs$.  Since $A_\Total=A_\Wobs$ the first condition, that $\alpha(\mathfrak{X}(A)_\Wobs,...,\mathfrak{X}(A)_\Wobs)\in A_\Wobs$ is void. The second condition states that $\alpha(X,Y_1,...,Y_k)\in A_0$ for any $X\in \mathfrak{X}(A)_\Null$, $Y_j\in \mathfrak{X}(A)_\Wobs$. By the previous Lemma we know that $X=\sum f_iX_i$ with $f_i\in A_0$, hence 
$$\alpha(X,Y_1,...,Y_k)=\sum_if_i\alpha(X_i,Y_1,...,Y_k)\in A_0.$$  
\end{proof}

For $\mathfrak{X}(A)_\Wobs$, the situation is a little more complicated.  For any vector field $X\in \X(M)$, the condition that the flow of $X$ preserves a closed subset $S$ implies $X(I_S)\subset I_S$, since $X(f)$ can be computed using the flow of $X$. The converse direction is not valid in general.
\begin{example}\label{ex:algtang} Let $M=\mathbb R$ and $S=[-1,1]$. Any function $f$ vanishing on $S$ also satisfies that $df|_S=0$, hence any vector field $X\in\mathfrak X(M)$ satisfies $X(I_S)\subset I_S$; however, the flow of a constant non-zero vector field will never preserve the set $S$.
\end{example}
\noindent
Crucially, preserving $S$ and preserving $I_S$ become equivalent when $S$ is a closed submanifold of $M$, cf. \autoref{lem:algeomtgt}.

Finally, we introduce two additional ingredients:
\begin{itemize}
	\item a Lie subalgebra $F \subset \mathfrak{X}(A)_\Wobs$, i.e., an involutive family of vector fields preserving the ideal $I_S$. The corresponding Lie--Rinehart algebra is $\mathcal{F} = \mathcal{F}_F = C^\infty(M) \cdot F + I_S \mathfrak{X}(M)$;
	\item a pre-multisymplectic form $\omega \in \Omega^{k+1}(M)$ such that $\omega\big(F, \mathfrak{X}(A)_\Wobs, \ldots, \mathfrak{X}(A)_\Wobs\big)\Big|_S = 0$.
\end{itemize}
Under these assumptions,  the main construction of \cite{blacker2023reduction} arises as a corollary of \autoref{thm:mainredgenalgleotired}.

\begin{corollary}[\cite{blacker2023reduction}]\label{cor:howtorecoverourreduction} In the above setting, the reduced $L_\infty$-algebra of observables is given by
	
	\begin{align}\label{eq:redsum}
		\frac{\Ham(B',\omega)_\Wobs}{\Ham(B',\omega)_\Null}=
		\frac{
			\left\lbrace
			\begin{array}{l}
				\alpha\in \Omega^{\leq k-2}(M) ~\mathrm{, respectively}\\	
				(\alpha,X) \in \Omega^{k{-}1}(M) \oplus\X(M)
			\end{array}
			\middle\vert
			\begin{array}{l l}
				\iota_\xi \alpha \in \CE(\mathfrak{X}(A))_\Null &\\
				\Lie_\xi\alpha \in \CE(\mathfrak{X}(A))_\Null&~\forall \xi \in F\\
				\hline
				\iota_X \omega = -\d \alpha &\\
				X \in \mathfrak{X}(A)_\Wobs &\\	
				 \left[ \xi, X \right] \in \mathcal{F} 
				 & \forall \xi \in F
			\end{array}
			\right\rbrace
		}{
			\left\lbrace
			\begin{array}{l}
				\alpha\in \Omega^{\leq k-2}(M) ~\mathrm{, respectively}\\	
				(\alpha,X) \in \Omega^{k{-}1}(M) \oplus\mathfrak{X}(M)
			\end{array}
			\middle\vert
			\begin{array}{l}
				\alpha \in \CE(\mathfrak{X}(A))_\Null\\
				\hline
				\iota_X \omega = -\d \alpha\\
				X \in \mathcal{F}
			\end{array}
			\right\rbrace
		}~.
	\end{align}
\end{corollary}

Let us illustrate this procedure with a few examples.

\begin{example}[covariant momentum map]\label{ex:covredalg}
Let us return to the setting of \autoref{ssec:symred}. 
I.e. assume we are given
\begin{itemize}
	\item $(M,\omega)$ a pre-multisymplectic manifold;
	\item $\mathfrak g\to \mathfrak X(M)$ an infinitesimal action preserving $\omega$;
	\item $\mu:\mathfrak g\to Ham^0(M,\omega)$ a covariant momentum map (\autoref{def:covmom}).
\end{itemize}

	We consider the closed set 
	$$
		S=
		\{x\in M ~\mid ~\mu_x(\xi)=0\in \Lambda^{k-1}T^*_xM~\forall \xi \in \mathfrak{g}\}
	$$
	and $F=\underline {\mathfrak{g}}$ the Lie algebra of infinitesimal generators of the action induced by $\mathfrak{g}$. 

	To carry out the reduction procedure from this section, we need to verify that  $\iota_\xi\omega \in \CE(\mathfrak{X}(A))_\Null$ for all $\xi\in F= \underline{\mathfrak{g}}$. 
	For this, we observe that by construction $\mu(\xi)$ vanishes on $S$, i.e. $\mu(\xi)\in \CE(\mathfrak{X}(A))_\Null$. But $\CE(\mathfrak{X}(A))_\Null$ is preserved by $d$, so  $\iota_\xi\omega=-d\mu(\xi)\in \CE(\mathfrak{X}(A))_\Null$. 
	\\
	Due to the nature of covariant momentum maps, the condition $\Lie_\xi\alpha\in \CE(\mathfrak{X}(A))_\Null$ can be omitted from \autoref{eq:redsum} for $(\alpha, X)\in \Ham^0(B'\omega)_\Wobs$. This follows from the equation:
	$$
		\Lie_\xi\alpha=\d\iota_\xi\alpha+\iota_\xi \d\alpha=\d\iota_\xi\alpha-\iota_\xi\iota_X\omega=\d\iota_\xi\alpha+\iota_X \d\mu(\xi)
	$$
	The statement now follows from $\alpha\in\CE(\mathfrak{X}(A))_\Null$, $\mu(\xi)\in\CE(\mathfrak{X}(A))_\Null$, and the fact that $\iota$ and $\d$ are constraint morphisms.
\end{example}

\begin{example}[Presymplectic case]
	In the symplectic case, certain parts of Equation \eqref{eq:redsum} are simplified automatically due to degree reasons. 
	The only relevant degree is (form degree) 0, hence the $\iota_\xi\alpha$ condition is void. Moreover, $CE(\mathfrak{X}(A))_\Null^0=I_S$, so we get:

\begin{align*}
		\frac{\Ham(B',\omega)_\Wobs}{\Ham(B',\omega)_\Null}=
	\frac{
		\left\lbrace
		\begin{array}{l}
			(\alpha,X) \in C^\infty(M) \oplus\X(M)
		\end{array}
		~\left\vert~
		\begin{array}{l l}
			\Lie_\xi\alpha \in I_S&~\forall \xi \in F\\
			\iota_X \omega = -\d \alpha &\\
			X (I_S)\subset I_S&\\
			 \left[ \xi, X \right] \in \mathcal{F} 
			& \forall \xi \in F
		\end{array}
		\right.
		\right\rbrace
	}{
		\left\lbrace
		\begin{array}{l}
			(\alpha,X) \in C^\infty(M) \oplus\X(M)
		\end{array}
		~\left\vert~
		\begin{array}{l}
			\alpha \in I_S\\
			\iota_X \omega = -\d \alpha\\
			X \in C^{\infty}(M)\cdot F + I_S\mathfrak X(M)
		\end{array}
		\right.
		\right\rbrace
	}~.
\end{align*}
\end{example}

To finish, let us have a look at what happens in the most classical symplectic case, combining the above two examples:

\begin{example}[Symplectic case with momentum]
	Let $(M,\omega)$ be a symplectic manifold and $\mathfrak g\to \mathfrak X(M)$ a symplectic action with (co)momentum map $\mu:\mathfrak g\to C^\infty(M)$. In this setting, a Hamiltonian pair $(\alpha, X)$ is completely determined by the function $\alpha \in C^\infty(M)$.
	\\
	Let us start by looking at a Hamiltonian pair in the denominator, i.e. $(\alpha, X)\in \Ham(B',\omega)_0$.
	The field $X=X_\alpha$ lies in $C^\infty(M)\underline{\mathfrak{g}}+I_S\mathfrak X(M)$, i.e. $X=\sum g_i\underline{\xi_i}+R$ with $R\in I_S\mathfrak X(M)$ and $g_i$ arbitrary smooth functions. 
	Contracting this with $\omega$ we get  
	$$\d\alpha=\sum -g_i \d\mu(\xi_i)-\iota_R\omega=\sum \d(-g_i\mu(\xi_i)) +(\d g_i)\wedge \mu(\xi_i)-\iota_R\omega~.
	$$
	The first summand is the exterior derivative of $\beta=-\sum (g_i\mu(\xi_i))$ in  $I_\mu$, where $I_\mu$ is the ideal generated by $\mu(\xi), \xi\in\mathfrak g$ (cf. \cite{SniatWein83}). 
	The second term lies in $I_S\,\Omega^1(M)$ as does the third, so $\d(\alpha-\beta)\in I_S\Omega^1(M)$. At the same time $\alpha-\beta$ lies in  $I_S\Omega^1(M)$, so it is a function which vanishes on $S$ along with its first derivative, i.e. it is an element of $Q=\bigcap_{x\in S}I_x^2$. This means that we can translate the conditions on $(\alpha,X_\alpha)$ into $\alpha\in \subset I_\mu+Q$. The same reasoning can be applied to $\{\xi,X\}$ in the numerator, hence we arrive at the expression:
	$$	
		\frac{\Ham(B',\omega)_\Wobs}{\Ham(B',\omega)_\Null}	=
		\frac{\{f\in C^{\infty}(M)~|~\{f,I_S\}\subset I_S, \{f,I_\mu\}\subset I_\mu+Q   \}}{I_\mu+Q}
	$$
Somewhat surprisingly, this construction (which a priori only yields a Lie algebra) ends up being a Poisson algebra. For a comparison of this reduction to other algebraic reduction schemes, we refer to our previous work \cite{blacker2023reduction}.
\end{example}

\subsection{Relation to geometric reduction and the residue defect}\label{ssec:geoRed-residue-defect}

In this subsection, we will study what happens to our algebraic constructions when geometric reduction \'a la Blacker~\cite{blackerReductionMultisymplecticManifolds2021} is possible. 
Before turning to the comparison between geometric and algebraic reduction, we will verify that the reductions of our constraint calculus recover classical Cartan calculus in these cases.

\subsubsection{Cartan Calculus for constraint manifolds}
\label{ss:geomcartan}
We now examine how the algebraic structures of Cartan calculus specialise in geometrically well-behaved contexts.
In particular, we will see that our Cartan calculus generalises the Cartan calculus for constraint manifolds presented in \cite{dippelkern}. Throughout this subsection, we will fix a constraint manifold, i.e. a triple $(M,S,D)$ where:
\begin{itemize}
	\item $M$ is a smooth manifold
	\item $S\subset M$ is a closed embedded submanifold
	\item $D\subset TS$ is a regular foliation, such that the quotient $S/D$ is again a smooth manifold.
\end{itemize}

We start by recalling that smooth functions/ sections can always be extended from $S$ to $M$:
\begin{lemma}\label{lem:extsec}Let $E\to M$ be a vector bundle, and $\phi\in \Gamma(E|_S)$. Then there exists a section $\Phi\in\Gamma(E)$ such that $\Phi|_S=\phi$.
\end{lemma}
\begin{proof}One way to prove it is by applying the tubular neighbourhood theorem (cf. e.g. \cite[Theorem 6.24]{leeIntroductionSmoothManifolds2012}) to get a local extension, and then using closedness of $S$ to globalise it (\cite[Lemma 10.12]{leeIntroductionSmoothManifolds2012}).
\end{proof}

By applying the above Lemma to the trivial $\mathbb{R}$-bundle, we immediately obtain that any function on $S$ can be extended to $M$. As a consequence:

\begin{corollary}\label{cor:Ared-constr-manifold}
	Let $A=(C^\infty(M), C^\infty(M), I_S)$ be the constraint algebra of smooth functions on a constraint manifold $(M, S, D)$. 
	Then the reduced algebra of observables is given by 
	$$
	A_\red=\frac{C^\infty(M)}{I_S}=C^\infty(S)
	$$
\end{corollary}

The second advantage of $S$ being an embedded submanifold is that a vector field $X\in\mathfrak{X}(M)$ being algebraically tangent to $S$ is equivalent to being geometrically tangent to $S$, i.e. the phenomenon described in \autoref{ex:algtang} can not arise.  

\begin{lemma}\label{lem:algeomtgt}
	The constraint algebra of derivations of $A$, as per \autoref{cor:Ared-constr-manifold}, satisfy:
	$\mathfrak{X}(A)_\Wobs=\{X\in\mathfrak X(M)~|~X|_S\in\mathfrak{X}(S)\}$. In particular,
	$$
	\mathfrak{X}(A)_\red=\mathfrak{X}(S)
	$$
\end{lemma}
\begin{proof} 
The implication $X|_S\in\mathfrak{X}(S)\implies X(I_S)\subset I_S$ is true in general. The converse direction can be shown by finding for any $s\in S, v\in T_SM\backslash T_sS$ a function $f\in I_S$ such that $\d f(v)\neq 0$. To find such a function, one can use the fact that locally $S$ looks like a hyperplane in $\mathbb {R} ^n$. The second part follows from the first and \autoref{lem:vanonsubset}.
\end{proof}

Together with \autoref{cor:bnbt} and \autoref{lem:extsec}, this implies that differential forms also behave as expected:

\begin{corollary} 
	Let us write $j$ for the inclusion $S\hookrightarrow M$. 
	Then we have:
 \begin{align*}
 	B_\Null=\CE(\mathfrak{X}(A))_\Null=\{\alpha\in\Omega(M)~|~j^*\alpha=0\}~.
 \end{align*}
Moreover, $	\CE(\mathfrak{X}(A))_\red=\Omega(S)$. 
\end{corollary}

Let us now look at the constructions involving the symmetries $D$. The constraint manifold $(M, S, D)$ models the quotient $M_\red=\frac{S}{\sim_D}$, i.e. the quotient of the closed subset $S$ by the foliation $D$, which we suppose to be a smooth manifold. We first define the Lie--Rinehart algebra $\mathcal F\subset\mathfrak{X}(M)$ from the distribution $D$ as follows:
\begin{align*}
	\mathcal F=\{X\in\mathfrak X(M)~|~X|_S\in\Gamma(D)\}~.
\end{align*}

We have that $\mathfrak{X}(A)_\Null=I_S\cdot\mathfrak X(M)\subset \mathcal F$. 
The space $A'_\Wobs$ can be written as 
\begin{align}\label{eq:normalsinconstrmfcase}
	A'_\Wobs=\{f\in C^\infty(M)~|~\d f(D)=0\}~.
\end{align}
This means that $A'_\Wobs$ contains exactly the functions that are $D$-invariant on $S$, i.e. which induce smooth functions on the quotient $M_\red$. Using \autoref{eq:normalsinconstrmfcase}, one can deduce that $\overline{\mathcal{F}}=\mathcal{F}$. Moreover, a direct computation in local coordinates shows that the normalizer of $\mathcal F$ in $\mathfrak{X}(A)_\Wobs$ contains exactly the vector fields whose restriction to $S$ is projectible to $M_\red$. Finally, $B'_\Wobs$ are exactly the basic forms. In total, we have: 

\begin{proposition} For constraint manifolds the spaces $\mathfrak Y_\red$ and $B'_\red$ are just the vector fields and forms on $M_\red$, i.e.
	\begin{align*}
		&\mathfrak{Y}_\red=\frac{\mathfrak{Y}_\Wobs}{\mathfrak{Y}_\Null}=\frac{\{X\in\mathfrak X(M)~|~X|_S\in \mathfrak{X}(S), [X|_S,\Gamma(D)]\subset \Gamma(D)\}}{\{X\in\mathfrak X(M)~|~X|_S=0\}}=\mathfrak X(M_\red)\\
		&B'_\red=\frac{B'_\Wobs}{B'_\Null}=\frac{\{\alpha\in\Omega(M)~|~\iota_{\Gamma(D)}\alpha|_S=0,  \Lie_{\Gamma(D)}\alpha|_S=0\}}{\{\alpha\in\Omega(M)~|~j^*\alpha=0\}}=\Omega(M_\red)\\
	\end{align*}
\end{proposition}

We close this subsection by relating the above with the Cartan calculus developed in \cite{dippelkern}. Our algebra of smooth functions on the constraint manifold corresponds to the one described in \cite[Proposition 3.5]{dippelkern}. The constraint derivations $\mathfrak{Y}=\mathfrak{X}(A)'$ then correspond to constraint vector fields by \cite[Proposition 3.45]{dippelkern} and our module $B'$ to their $C\Omega_\boxtimes^\bullet$. The identities of the constraint Cartan calculus given by the BV-module structure, then reproduce \cite[Proposition 3.48]{dippelkern}. 

\begin{remark} In \cite{dippelkern} two versions of the space of differential forms are considered, one based on the weak tensor product $\otimes$ and the other on the strong tensor product $\boxtimes$. They show that only the latter yields a well-behaved Cartan calculus. In particular, the weak version fails to support a consistent de~Rham differential (see~\cite[Prop.~3.47, Ex.~3.48]{dippelkern}). For the same reason, we chose weak tensor product for the exterior algebra of a Lie--Rinehart algebra, leading to a strong constraint algebra structure on the dual, cf. \autoref{ex:constextalgdual}.
\end{remark}

\subsubsection{Multisymplectic reduction with respect to a covariant momentum map}

Let $(M,\omega)$ be a multisymplectic manifold with the action of a connected Lie group $G$, and $\mu$ a covariant moment map for this action (see \autoref{def:covmom}), such that its vanishing locus $S$ is an embedded submanifold and the action of $G$ on $S$ is free and proper. 
Then we can apply \cite[Theorem 1]{blackerReductionMultisymplecticManifolds2021}, to obtain a pre-mulisymplectic manifold $(M_{\red},\omega_{\red})$, where $M_{\red}={S}\slash {G}$ and $\omega_{\red}$ is compatible with $\omega$ when pulled back to $S$. We note that, since $G$ is connected, its fundamental vector fields form a regular foliation when restricted to $S$; in particular, we are in the setting of the previous subsection.

In this setting, we have two observable $L_\infty$-algebras we could consider for the quotient, the geometric one $\Ham(M_{\red},\omega_{\red})$ and the algebraic reduction described in \autoref{ex:covredalg}.
If a vector field $X$ and a form $\alpha$ satisfy $\iota_X\omega=-\d \alpha$ and both pass to the quotient, then their reductions satisfy the same relation with respect to $\omega_{\red}$, implying that there is a natural injective morphism 
\begin{align}\label{eq:redcomp}
	\frac{\Ham(B',\omega)_\Wobs}{\Ham(B',\omega)_\Null}\to \Ham(M_{\red},\omega_{\red}).
\end{align}
This morphism preserves the $L_\infty$-structures. It is, however, in general not surjective.

\begin{example}
	Let us consider $M=\mathbb R^5$ with the coordinates $t,x_1,x_2,y_1,y_2$ and the multisymplectic form $\omega=\d t\wedge (\d x_1\wedge \d y_1+\d x_2\wedge \d y_2)$. We consider the 1-dimensional action by $\underline{\xi}=\partial_{x_1}$ and the momentum $\mu(\xi)=-y_1\d t$. We get $S=\{y_1=0\}$ and $M_{\red}=\mathbb R^3$ with coordinates $t,x_2,y_2$ and multisymplectic form  $\d t\wedge \d x_2\wedge \d y_2$.
	We can now consider the reduction of the Hamiltonian pair $(x_2 t \d y_2,x_2\partial_t-t\partial_{x_2})$. It can never be the (multisymplectic observable) reduction of any Hamiltonian pair on $M$. The vector field $x_2\partial_t-t\partial_{x_2}$ can be seen as a reducible vector field on $M$; any reducible vector field with the same reduction will have the form:
	$$
	V=x_2\partial_t-t\partial_{x_2}+ f\partial_{x_1}+ y_1\cdot R
	$$
	for an unknown function $f$ and vector field $R$. 
	However, $\iota_V\omega$ can never be closed for such a vector field, because $\d\iota_V\omega$ would always have a non-trivial $\d x_2\wedge \d x_1  \wedge \d y_1$-component, over $y_1=0$.
\end{example}

The failure of the morphism \eqref{eq:redcomp} to be surjective is called \emph{residue defect} (\cite{blacker2023reduction}) and can be explained as follows: In our setup we are demanding a very strict condition on Hamiltonian pairs, we demand $\d\alpha = -\iota_X\omega$ globally on $M$, not only on $S$, (and without relaxing the equation in $F$-directions). 
In the symplectic (Marsden-Weinstein-Meyer) case, this is not a big deal, since the form uniquely determines the vector field both on $M$ and on $M_{\red}$. However, in the multisymplectic case, Hamiltonian pairs might be scarce on $M$ and become less rare on $M_{\red}$, which is precisely what happens in the above example. 
Note that the non-surjectivity only happens on the level of $\Ham^0$, where we have pairs of vector fields and forms, rather than just forms in the structure.
\\
In order to circumvent this problem, we could look at Hamiltonian pairs $(\alpha, X)$, where $\d\alpha+\iota_X\omega\in B'_0$; however, such couples are not elements in the initial $L_\infty$-algebra on $M$, thus should not be seen as elements in the reduction. \\

In fact, we could generalise even further: We could go through the whole $L_\infty$-algebra from a cocycle construction in the general constraint setting: We would start, e.g. with a constraint cocycle ($\omega\in B'_\Wobs$ such that $\d\omega\in B'_\Null$) and construct operations starting from it. However, in this setting, we would have an $L_\infty$-algebra only on the quotient; the structure on the $\Total$, $\Wobs$, and $\Null$ components would be more general. This more general structure could be an interesting subject of study in its own right.


			\bibliographystyle{hep}
			\bibliography{biblio}

\end{document}